\pdfoutput=1

\documentclass[letterpaper, 12pt, reqno]{amsart}
\usepackage{amsmath, amsfonts, amssymb, amsthm, array, stmaryrd, graphicx, hyperref, mathrsfs, eucal, caption, subcaption, soul, cancel, ulem, float}

\usepackage[alphabetic]{amsrefs}
\AtBeginDocument{%
	\def\MR#1{}
}

\usepackage{mathptmx}
\DeclareSymbolFont{Symbols}{OMS}{zplm}{m}{n}
\DeclareMathSymbol{\infty}{\mathord}{Symbols}{"31}

\usepackage[hhmmss]{datetime}

\allowdisplaybreaks

\usepackage[letterpaper, total={6.5in, 8.5in}, left=1in, top=1.25in]{geometry}
\newcount\Comments  
\usepackage{color}
\definecolor{burntorange}{RGB}{204, 85, 0}
\definecolor{oregongreen}{RGB}{0, 112, 48}
\definecolor{oaklandgold}{RGB}{181, 154, 87}
\definecolor{sydneyblue}{RGB}{1, 72, 164}

\newcommand{\kibitz}[2]{\ifnum\Comments=1\textcolor{#1}{#2}\fi}

\theoremstyle{plain}

\newtheorem{lemma}{Lemma}[section]

\newtheorem{conj}[lemma]{Conjecture}
\theoremstyle{definition}

\theoremstyle{remark}
\newtheorem{remark}[lemma]{Remark}

\numberwithin{equation}{section}

\newcommand{\p}{\partial}

\newcommand{\bn}{\begin{enumerate}}
	\newcommand{\en}{\end{enumerate}}
\newcommand{\bi}{\begin{itemize}}
	\newcommand{\ei}{\end{itemize}}
\newcommand{\bqq}{\begin{eqnarray*}}
	\newcommand{\eqq}{\end{eqnarray*}}
\newcommand{\balg}{\begin{align*}}
\newcommand{\ealg}{\end{align*}}

\DeclareMathOperator{\Rm}{Rm}

\DeclareMathOperator{\cyl}{Cyl}

\begin{document}
\title[Numerical stability analysis of noncompact Type-II MCF solutions]{A numerical stability analysis of mean curvature flow of noncompact hypersurfaces with Type-II curvature blowup}

\author{David Garfinkle}
\address{Department of Physics, Oakland University, Rochester, MI 48309, USA}
\email{garfinkl@oakland.edu}

\author{James Isenberg}
\address{Department of Mathematics, University of Oregon, Eugene, OR 97403, USA}
\email{isenberg@uoregon.edu}

\author{Dan Knopf}
\address{Department of Mathematics, The University of Texas at Austin, Austin, TX 78712, USA}
\email{danknopf@math.utexas.edu}

\author{Haotian Wu}
\address{School of Mathematics and Statistics, The University of Sydney, NSW 2006, Australia}
\email{haotian.wu@sydney.edu.au}


\keywords{Mean curvature flow; Type-II singularities; noncompact hypersurfaces; stability analysis; numerical method}

\subjclass[2010]{53C44, 35K59, 65M06, 65D18}

\begin{abstract}
We present a numerical study of the local stability of mean curvature flow of rotationally symmetric, complete noncompact hypersurfaces with Type-II curvature blowup. Our numerical analysis employs a novel overlap method that constructs ``numerically global'' (i.e., with spatial domain arbitrarily large but finite) flow solutions with initial data covering analytically distinct regions. Our numerical results show that for certain prescribed families of perturbations, there are two classes of initial data that lead to distinct behaviors under mean curvature flow. Firstly, there is a ``near'' class of initial data which lead to the same singular behaviour as an unperturbed solution; in particular, the curvature at the tip of the hypersurface blows up at a Type-II rate no slower than $(T-t)^{-1}$. Secondly, there is a ``far'' class of initial data which lead to solutions developing a local Type-I nondegenerate neckpinch under mean curvature flow. These numerical findings further suggest the existence of a ``critical'' class of initial data which conjecturally lead to mean curvature flow of noncompact hypersurfaces forming local Type-II degenerate neckpinches with the highest curvature blowup rate strictly slower than $(T-t)^{-1}$.
\end{abstract}

\maketitle


\section{Introduction}\label{intro}

The nature of finite time singularities in geometric flows such as Ricci flow and mean curvature flow (MCF) has long been a focus of interest in the study of these flows. In both flows, the formation of Type-I singularities, which are characterized by the supremum of the product of the finite time to the singularity and the curvature norm being finite, predominate. Much less prevalent are Type-II singularities, in which case the supremum of the product of the time to the singularity and the curvature norm is necessarily infinite.

For immersed one-dimensional MCF (i.e., curve shortening flow), Type-II behavior can occur for solutions originating from open sets of initial data \cite{EW87}. Nonetheless, it is a folklore conjecture that Type-I behavior is generic for compact (embedded) solutions. For MCF, there is in fact strong evidence in important work of Colding and Minicozzi \cite{CM12} that only a restricted subclass of Type-I tangent flows (self-shrinkers) is generic. For example, they prove that if $\Sigma$ is a smooth, complete, embedded self-shrinking surface without boundary and with polynomial volume growth that is \emph{not} a generalized cylinder, then there is a graph $\tilde\Sigma$ over $\Sigma$ of a compactly supported function with arbitrarily small $C^m$-norm (for any fixed $m$) such that $\Sigma$ cannot occur as the tangent flow of MCF originating from $\tilde\Sigma$.

\begin{figure}[H]
	\includegraphics[width=0.6\textwidth]{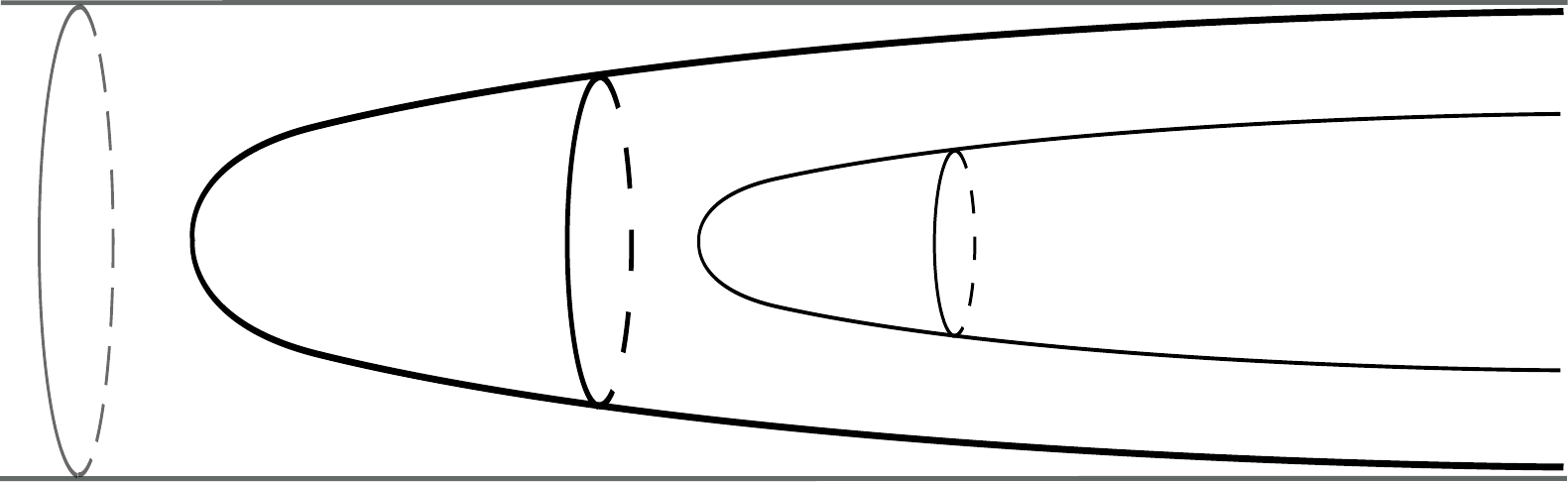}
	\caption{Diagram of a MCF solution forming a degenerate neckpinch at spatial infinity with the curvature at the tip (on the left) blowing up at a Type-II rate, as constructed in \cite{IW19,IWZ20}.}\label{fig1}
\end{figure}

In the noncompact setting, Type-II singularities can form in solutions originating from an open set within the class of rotationally symmetric initial data for embedded MCF, as is shown by two of the authors of this paper and their collaborator \cite{IW19, IWZ20} (similar behavior in Ricci flow has been verified in \cite{Wu14}). They prove that MCF of embedded noncompact hypersurfaces satisfying certain conditions\footnote{Precisely, the hypersurface is a complete, rotationally symmetric, strictly convex graph over a shrinking ball that is asymptotically enveloped within a shrinking cylinder near spatial infinity; see Figure \ref{fig1}.} necessarily develop a Type-II singularity. Furthermore, the maximum of the curvature measured in terms of $|h|$, where $h$ denotes the second fundamental form, of such a solution occurs at the tip, with the Type-II blowup of the flow approaching a translating solution known as the bowl soliton --- the unique (up to rigid motion) translating solution that is rotationally symmetric and strictly convex \cite{Has15} --- in the region surrounding the tip (i.e., the left most point on the hypersurface in Figure \ref{fig1}). At spatial infinity (i.e., on the far right of Figure \ref{fig1}), the solution remains asymptotic to a shrinking cylinder and hence forms a Type-I singularity. In particular, such a solution is an example of a degenerate neckpinch in MCF forming at spatial infinity.

What is the stability of the noncompact mean curvature flow behavior illustrated in Figure \ref{fig1} and discussed in \cite{IW19} and \cite{IWZ20}? In this work, we carry out a numerical stability analysis by simulating mean curvature flows of embedded noncompact hypersurfaces that are perturbations of the embeddings considered in \cite {IW19, IWZ20}. Retaining the rotational symmetry in these perturbations, we find that there is one class of such perturbations for which the simulated mean curvature flows develop Type-II singularities modelled by the bowl soliton, and there is another set of these perturbations for which the numerical simulations develop Type-I singularities. For this second class, the Type-I singularities that occur are all of the ``neckpinch'' type. The perturbations that lead to mean curvature flows with Type-II singularities generally involve distortions of the initial embeddings very close to the tip (which we call the ``near class'' in this paper), while those perturbations that lead to Type-I singularities generally involve distortions that are not close to the tip (called the ``far class'' in this work). Interestingly, as noted below in Figure \ref{near4}, there are perturbations of the initial data which are not especially close to the tip, yet the numerical simulation of the MCF starting from that perturbed initial data set approaches the bowl soliton with Type-II behavior at the tip and with no neckpinch singularity forming away from the tip. We define the ``near'' and ``far'' classes (in Section \ref{near_class} and \ref{far_class}, resp.) based on the MCF behavior rather than on the location of the perturbation.

One expected consequence of the aforementioned numerical results is that if one were to consider a one-parameter family (with parameter $s$) of rotationally symmetric perturbed initial embeddings that proceeds from the near class to the far class, then the MCF solutions evolving from all the embeddings up to a certain value $s_0$ develop Type-II singularities, while the flows evolving from the embeddings up to a certain value $s_1\geq s_0$ develop Type-I neckpinch singularities. The question then arises regarding what happens for the MCF solution that evolves from an embedding corresponding to an intermediate parameter value $s\in [s_0,s_1]$.

For analogous numerical simulation studies of Ricci flow on a compact manifold carried out by two of the authors \cite{GI05, GI08}, one set of parameters lead to nondegenerate neckpinches while another set of parameters lead to round singularities, and a very distinctive behavior --- a degenerate neckpinch, corresponding to some intermediate, critical parameter value $s_0$ --- is observed. The same phenomenon  is seen in the remarkable numerical studies of gravitational collapse carried out by Choptuik. Notably, in certain of his numerical simulations of the gravitational collapse governed by solutions of Einstein’s gravitational field equations, very distinctive critical behavior is found to occur for a unique choice of the initial data along a one-parameter family of such gravitational data \cite{Chop93}. Interestingly, for an alternative collection of Choptuik’s numerical simulations of gravitational collapse, critical behavior is not found \cite{Gund-MG07}.

One challenge in carrying out these numerical simulations is the need to ensure compatibility of the evolutions that we carry out in two different regions with two different coordinate systems. Using the coordinate $z$ for the direction parallel to the enveloping cylinder, and the coordinate $r$ for the radial direction orthogonal to $z$, we are led to evolve the function $z(r,t)$ in the region near the tip, and to evolve the function $r(z,t)$ in the region away from the tip and along the cylindrical region. Coordinate transformations of this sort are of course familiar in differential geometry, but they do present a challenge in these numerical simulations. Since the near class of initial data generally involves perturbations in the region where we work with $z(r,t)$ while the far class of initial data generally involves perturbations in the region where we work with $r(z,t)$, we note that these coordinate transformations make it somewhat difficult to very closely investigate whether distinctive critical behavior occurs for perturbations located at the transition region.

While we have looked very carefully at numerical simulations of the mean curvature flows for initial embeddings very close to the transition from the near class to the far class, we have not found any distinctive critical behavior there. This may mean that the potential critical behavior does not exist, or that it does exist but is too unstable to be readily found by numerical experiments. Inferring from our numerical results, we conjecture that a critical behavior exists (see Conjecture \ref{conj_critical}), and further investigations are needed.

The perturbations of the initial embeddings we consider in this work are all rotationally symmetric. In future work, we plan to consider perturbations that do not retain this symmetry. We conjecture that if these perturbed initial embeddings remain small and confined to a region near the tip, the resulting MCF solutions are likely to still produce Type-II singularities modeled by a bowl soliton in the neighborhood of the tip.

This paper is organized as follows. In Section 2, we briefly review from \cite{IW19, IWZ20} the construction of noncompact, rotationally symmetric MCF solutions with Type-II curvature blowup at the tip. In Section \ref{num_method}, we discuss the computational techniques used to carry out the numerical simulations done here. Numerical results for the near classes and the far classes are presented in Section \ref{num_results}. These results lead to conjectures on the asymptotic behavior of mean curvature flows originating from initial data similar to the near class, and the existence of a critical class as one interpolates between the near classes and the far classes. Concluding remarks and discussions of future directions appear in Section \ref{discussion}.

\section*{Acknowledgements}

DG thanks the National Science Foundation for support in PHY-1806219. JI thanks the National Science Foundation for support in grant PHY-1707427. DK thanks the Simons Foundation for support in Award 635293. HW thanks the Australian Research Council for support in DE180101348.

\section{Set up}\label{setup}

Following the set up in \cite{IW19, IWZ20}, we consider the mean curvature flow of rotationally symmetric hypersurfaces embedded in Euclidean space. For any point $(x_0,x_1,\ldots,x_n)\in\mathbb{R}^{n+1}$ for $n\geqslant 2$, we write
\begin{align*}
z = x_0, \quad\quad  r = \sqrt{x_1^2 + \cdots + x_n^2}.
\end{align*}
We consider a noncompact hypersurface $\Gamma$ which is obtained by rotating the graph of $r(z)$, $a\leq z<\infty$, around the $z$-axis. We assume that $r(z)$ is strictly concave so that $\Gamma$ is strictly convex, and that $r$ is strictly increasing with $r(a)=0$ and with $\lim\limits_{z\nearrow \infty}r(z) = r_0$, where $r_0$ is the radius of the enveloping cylinder. The function $r$ is assumed to be smooth except at $z=a$. We note that this particular non-smoothness of $r$ is a consequence of the choice of the (cylindrical-type) coordinates; if the time-dependent flow function $r(z,t)$ is inverted in a particular way, then this irregularity is removed. The point where $r=0$ is called the \emph{tip} of the hypersurface.

We focus our attention on the class $\mathscr{G}$ of complete hypersurfaces that are rotationally symmetric, (strictly) convex\footnote{Throughout this paper, ``convex'' means ``strictly convex''.}, smooth graphs over a ball and asymptotic to a cylinder. One readily verifies that embeddings with these properties are preserved by MCF (see for example \cite{SS14}). By the general existence result in \cite{SS14}, a MCF solution starting from any hypersurface in this class moves towards its open end (e.g., to the right in Figure \ref{fig1}) and remains asymptotic to a shrinking cylinder; moreover, the hypersurface disappears at spatial infinity at the same time as the cylinder collapses. We call this finite time the ``vanishing time'' of the hypersurface and denote it by $T$. The results of \cite{IW19} and \cite{IWZ20} exhibit, for each real number $\gamma\geq 1/2$, MCF
solutions in the class $\mathscr{G}$ with the Type-II curvature blowup rate $(T-t)^{-(\gamma+1/2)}$ at the tip, and also determine the precise geometry of the solution both near the tip and near spatial infinity.

Let $\Gamma_t$ denote a MCF solution constructed in \cite{IW19, IWZ20}. Representing $\Gamma_t$ by the profile of rotation, i.e., the graph of $r(z,t)$, one finds that the function $r$ satisfies the PDE
\begin{align}\label{r(z,t)}
r_t & = \frac{r_{zz}}{1+r^2_z} - \frac{n-1}{r}.
\end{align}
A multi-scale blowup analysis of equation \eqref{r(z,t)} is carried out in \cite{IW19} for the case $\gamma>1/2$, and in \cite{IWZ20} for the ``limiting'' case $\gamma=1/2$. We recall the relevant details for each case below.

\subsection{The case $\gamma>1/2$}
In the rescaled time and space parameters
\begin{align*}
\tau = -\log(T-t), \quad y = z(T-t)^{\gamma-1/2}, \quad \phi(y,\tau) = r(x,t)(T-t)^{-1/2},
\end{align*}
equation \eqref{r(z,t)} becomes
\begin{align}\label{phi(y,tau)}
\left. \p_\tau \right\vert_y \phi & = \frac{e^{-2\gamma\tau} \phi_{yy}}{1+e^{-2\gamma\tau}\phi^2_y} - (1/2-\gamma) y \phi_y
- \frac{(n-1)}{\phi} + \frac{\phi}{2},
\end{align}
where the notation $\left. \p_\tau \right\vert_y $ means taking the partial derivative in $\tau$ while keeping $y$ fixed. We readily note that equation \eqref{phi(y,tau)} admits the constant solution $\phi \equiv \sqrt{2(n-1)}$, which  corresponds to the collapsing cylinder (which is a self-shrinking solution of MCF).

We now introduce two more coordinate transformations. Firstly, because the hypersurface under consideration is a convex graph over a ball, it is useful to invert the coordinates and work with 
\begin{align*}
y( \phi,\tau) & = y\left( \phi(y,\tau), \tau \right).
\end{align*}
In terms of $y(\phi,\tau)$, the MCF equation, which is equivalent to equations \eqref{phi(y,tau)} and \eqref{r(z,t)}, becomes
\begin{align}\label{y(phi,tau)}
\left. \p_\tau \right\vert_{\phi} y & = \frac{y_{\phi\phi}}{1 + e^{2\gamma\tau} y^2_\phi} + \left(\frac{(n-1)}{\phi} - \frac{\phi}{2} \right) y_\phi + (1/2-\gamma) y.
\end{align}
Secondly, it is further useful to work with the quantity $\lambda := -1/y$, since by using $\lambda$, the asymptotically cylindrical end of the embedded hypersurface corresponding to large values of $y$ is effectively compactified. The evolution equation for $\lambda(\phi,\tau)$ is
\begin{align}\label{lambda}
\left. \p_\tau \right\vert_\phi \lambda & = \frac{\lambda_{\phi\phi}-2\lambda^2_\phi/\lambda}{1+e^{2\gamma\tau}\lambda^2_\phi/\lambda^4} + \left(\frac{n-1}{\phi} - \frac{\phi}{2}\right)\lambda_\phi + \left(\gamma-\frac{1}{2}\right)\lambda.
\end{align}
By rotational symmetry, we can let $\phi\in(-\sqrt{2n-2},\sqrt{2n-2})$ for $n\geq 2$, and $\tau\in[\tau_0,\infty)$ for some large $\tau_0=-\log(T-t_0)$ for the initial time $t_0$. The boundary conditions for \eqref{lambda} are $\lambda(-\sqrt{2n-2},\tau)=\lambda(\sqrt{2n-2},\tau)=0$ for all $\tau\geq\tau_0$.

Using the method of matched asymptotics, formal solutions to equation \eqref{lambda} can be derived in two regions: the interior region where $\zeta:=\phi e^{\gamma\tau} = \mathcal O(1)$, and the exterior region, which is the complement of the interior region (see details in \cite[Section 2]{IW19}). Global
initial data for \eqref{lambda} are then defined by joining the formal functions in the exterior and the interior regions. Precisely, taking any real number $c>0$ (e.g., $c=1$) and letting $A:= c(2n-2)^{\gamma-1/2}$, we define
\begin{align}\label{id}
\hat\lambda_0(\phi) := \left\{
\begin{array}{lr}
\begin{array}{l} -A + e^{-2\gamma\tau_0}F(\zeta)- e^{-2\gamma\tau_0}F(R_1) \\
\quad + \left[ A -c\left(2n-2-(R_1e^{-\gamma\tau_0})^2\right)^{\gamma-1/2} \right] 
\end{array}, & 0\leq|\zeta|\leq R_1,\\ \\
-c(2n-2-\phi^2)^{\gamma-1/2}, & R_1e^{-\gamma\tau_0}\leq|\phi|<\sqrt{2(n-1)},
\end{array}
\right.
\end{align}
where $R_1$ is some large constant, and $F$ is the unique solution to the ODE initial value problem
\begin{align*}
\frac{F_{\zeta\zeta}}{1+F^2_{\zeta}/A^4} + (n-1) F_{\zeta}/\zeta = (\gamma-1/2)A,\quad\quad F(0)=F_{\zeta}(0)=0.
\end{align*}
In fact, $F(\zeta)$ is defined for all $\zeta\geq 0$ and gives the profile function for a scaled copy of the bowl soliton. For any $R_1\gg1$ we choose, we can always find a sufficiently large $\tau_0$ so that $\phi = \zeta e^{-\gamma\tau_0}$ is close enough to $0$ for
definition \eqref{id} to make sense and yield piecewise-smooth continuous initial data for the PDE \eqref{lambda}. The unperturbed initial data sets are specified by the constants $(n, c, R_1, \tau_0)$. In the end of Section \ref{num_method}, we discuss how to specify the perturbed initial data sets used in our simulations.

\begin{remark}
	For the numerical study in this paper, which is carried out using standard finite differencing methods, piecewise-smooth continuous initial data for \eqref{lambda} are sufficient. An extra step to smooth the initial data is needed in \cite{IW19}.
\end{remark}

\subsection{The case $\gamma=1/2$}
For $a>0$ to be chosen, we consider the rescaled time and space parameters
\begin{align*}
\tau  &= -\log(T-t),\quad y = z + a \log(T-t),\quad \phi(y,\tau) = r(x,t)(T-t)^{-1/2}.
\end{align*}
Then equation \eqref{r(z,t)} becomes
\begin{align}\label{phi(y,tau)-crit}
\left. \p_\tau \right\vert_y \phi & = \frac{e^{-\tau} \phi_{yy}}{1+e^{-\tau}\phi^2_y} + a \phi_y
+ \frac{\phi}{2} - \frac{(n-1)}{\phi}.
\end{align}
Letting
\begin{align*}
y(\phi,\tau) & = y\left( \phi(y,\tau), \tau \right),\quad\quad \lambda(\phi,\tau):=-1/y(\phi,\tau).
\end{align*}
as in the case that $\gamma>1/2$, we find that $\lambda(\phi,\tau)$ evolves according to the PDE
\begin{align}\label{lambda-crit}
\left. \p_\tau \right\vert_\phi \lambda & = \frac{\lambda_{\phi\phi}-2\lambda^2_\phi/\lambda}{1+e^{\tau}\lambda^2_\phi/\lambda^4} + \left(\frac{n-1}{\phi} - \frac{\phi}{2}\right)\lambda_\phi - a \lambda^2.
\end{align}
By rotational symmetry, we consider $\phi\in(-\sqrt{2n-2},\sqrt{2n-2})$ for $n\geq 2$, and $\tau\in[\tau_0,\infty)$ for some large $\tau_0$. The boundary conditions are $\lambda(-\sqrt{2n-2},\tau)=\lambda(\sqrt{2n-2},\tau)=0$ for all $\tau\geq\tau_0$.

As in the case that $\gamma>1/2$, the formal solutions to equation \eqref{lambda-crit} are derived in two regions: the interior region where $\zeta:=\phi e^{\tau/2} = \mathcal O(1)$, and the exterior region, which is the complement of the interior region (see details in \cite[Section 2]{IWZ20}). Global initial data for~\eqref{lambda} are then defined by joining the formal functions in the exterior and the interior regions. Precisely, 
we take any $a>0$, fix any $c>a\log(2n-2)$, let $A:=1/\left(c-a\log(2n-2)\right)$, and define
\begin{align}\label{id-crit}
\widehat\lambda_0(\phi) := \left\{
\begin{array}{lr}
-A + e^{-\tau_0}F(\zeta)- e^{-\tau_0}F(R_1)\\ \quad +A-\left(c-a\log(2n-2-R_1^2e^{-\tau_0})\right)^{-1}, & 0\leqslant |\zeta|\leqslant R_1, \\ \\
-1/\left(c-a\log(2n-2-\phi^2)\right), & R_1e^{-\tau_0/2}\leqslant |\phi|<\sqrt{2(n-1)},
\end{array}
\right.
\end{align}
where $R_1$ is a large constant, and $F$ is the unique solution to the ODE initial value problem
\begin{align*}
\frac{F_{\zeta\zeta}}{1+F^2_\zeta/A^4} + (n-1) \frac{F_\zeta}{\zeta} = aA^2,\quad\quad F(0)=F_\zeta(0)=0.
\end{align*}
As remarked above, $F(\zeta)$ is the profile function for a scaled copy of the bowl soliton. For any $R_1\gg1$ we choose, we can find $\tau_0$ large enough so that $\phi = \zeta e^{-\tau_0/2}$ is close enough to $0$ for
definition \eqref{id-crit} to make sense and yield piecewise-smooth continuous initial data for the PDE~\eqref{lambda-crit}. The unperturbed initial data sets are specified by the constants $(n, c, R_1, \tau_0)$. In the end of Section \ref{num_method}, we discuss how to specify the perturbed initial data sets used in our simulations.

\begin{remark}
	As noted above, since the numerical study in this paper is carried out using finite differencing methods, the use of piecewise-smooth continuous initial data for \eqref{lambda-crit} is sufficient.
\end{remark}

\section{Numerical method}\label{num_method}

Our numerical method essentially consists of writing the mean curvature flow PDE in a standard parabolic form and then using a standard finite-difference method for parabolic equations.  However, for reasons to be explained below, we find it convenient to write the equation in two different ways and to use both of them.

We recall that equation \eqref{r(z,t)} tells us that, under mean curvature flow, the profile function $r(z,t)$ of a rotationally symmetric $n$-dimensional hypersurface embedded in $\mathbb{R}^{n+1}$ satisfies the following evolution equation:
\begin{align}\label{eq-r}
{\partial _t} {|_z} (r) & = \frac{r_{zz}}{1+r^2_z} - \frac{n-1}{r}.
\end{align}
If $r(z,t)$ is invertible, we can write $z(r,t)=z(r(z,t),t)$ and then $z(r,t)$ satisfies the PDE
\begin{align}\label{eq-z}
\left.\partial_t\right\vert_r(z) & = \frac{z_{rr}}{1+z^2_r} + (n-1) \frac{z_r}{r}.
\end{align}

In the rest of this paper, we consider the mean curvature flow of noncompact, rotationally symmetric surfaces in $\mathbb{R}^3$. (We briefly consider higher dimensions in Section \ref{n=3} and verify that those solutions display the same qualitative behaviors as we observe for the MCF for surfaces.) Since all of our simulations involve (hyper-)surfaces with rotational symmetry, the evolution equations \eqref{eq-r} and \eqref{eq-z} for higher dimensions differ from those for surfaces only in terms of a coefficient appearing on the right hand side; consequently the same numerical analysis applies. Consider a rotationally symmetric map from the two-sphere with coordinates $(\theta, \varphi)$ to
$\mathbb R^3$ with cylindrical coordinates $(r,z,\varphi)$; that is, $r$ and $z$ are functions of 
$\theta$.  This is the setup for mean curvature flow of a compact surface.  For the noncompact surfaces that are the main focus of this paper, we can consider maps from an open portion of the two-sphere.  In terms of Cartesian coordinates, we have
\begin{align*}
(x,y,z) = (r \cos \varphi, r \sin \varphi, z)
\end{align*}
Under mean curvature flow, the profile $r(z,t)$ evolves by (cf. \eqref{r(z,t)}) 
\begin{align}
{\partial _t} {|_z} (r) = {\frac {r_{zz}} {1 + {r_z ^2}}} \; - \; {\frac 1 r}. \label{mcf1}
\end{align}
Correspondingly, if we invert the coordinates and work with
\begin{align*}
z(r,t)= z(r(z,t),t),
\end{align*}
then under MCF, $z(r,t)$ satisfies the PDE
\begin{align}
{\partial _t} {|_r} (z) = {\frac {z_{rr}} {1 + {z_r ^2}}} \; + \; {\frac {z_r} r}. \label{mcf2}
\end{align}
Thus it appears that we have a choice of evolving either equation \eqref{mcf1} or equation \eqref{mcf2}.

Let us consider first what is involved in numerically evolving (\ref{mcf1}).  We approximate the 
function $r(t,z)$ by the values $r^k _i$ that the function takes at points $z_i$ equally spaced with spacing $\Delta z$ and with times $t_k$ equally spaced with spacing $\Delta t$.  We use standard centered finite differences, so that $r_z$ and $r_{zz}$ are approximated by
\begin{equation}
{r_z} = {\frac {{r^k _{i+1}}-{r^k _{i-1}}} {2 \Delta z}}
\label{fd1}
\end{equation}
and
\begin{equation}
{r_{zz}} = {\frac {{r^k _{i+1}} + {r^k _{i-1}} - 2 {r^k _i}} {{(\Delta z)}^2}},
\label{fd2}
\end{equation}
respectively.

Time evolution is carried out using the Euler method, so that 
\begin{equation}
{r^{k+1} _i} = {r^k _i} + \Delta t {\partial _t} r, 
\label{Euler1}
\end{equation}
where ${\partial _t} r$ is the finite difference version of the right hand side of equation (\ref{mcf1}) with the spatial derivatives evaluated using equations (\ref{fd1}) and (\ref{fd2}).  The standard von Neumann stability analysis of equation (\ref{Euler1}) reveals that the time step must satisfy the Courant condition
\begin{equation}
\Delta t < {\textstyle {\frac 1 2}} {{(\Delta z)}^2}.
\label{CFL1}
\end{equation}
Equation (\ref{CFL1}) is already a severe constraint on the type of numerical simulations that we can do. We need a small spatial step $\Delta z$ to accurately model the small spatial scales involved in the approach to the singularity.  However, the Courant condition then requires a very small time step.  The advantage of the rotational symmetry means that even with small $\Delta z$, only a moderate amount of computer memory is needed to store the values of $r$ at the spatial grid points, and only a moderate amount of computer time is needed to take the enormous number of time steps required by the tiny time step $\Delta t$.  

A more serious difficulty comes from the nature of the finite difference approximations in equations (\ref{fd1}) and (\ref{fd2}).  These approximations are good for functions that vary very little on the scale of $\Delta z$.  However, for our surfaces, $r_z$ and $r_{zz}$ are singular at the tip, which means that no matter how small we make $\Delta z$, the finite difference approximations fail at the tip.  

Alternatively, we can evolve equation (\ref{mcf2}) which is regular at the tip.  Here, the function $z(t,r)$ is represented by the values $z^k _i$ that it takes on a spatial grid with spacing $\Delta r$ and at times with spacing $\Delta t$. In this setup, the finite difference approximations used are
\begin{eqnarray}
{z_r} = {\frac {{z^k _{i+1}}-{z^k _{i-1}}} {2 \Delta r}},
\label{fd3}
\\
{z_{rr}} = {\frac {{z^k _{i+1}} + {z^k _{i-1}} - 2 {z^k _i}} {{(\Delta r)}^2}}.
\label{fd4}
\end{eqnarray}
Time evolution is again carried out using the Euler method, so that
\begin{equation}
{z^{k+1} _i} = {z^k _i} + \Delta t {\partial _t} z,
\label{Euler2}
\end{equation}
where ${\partial _t} z$ is the finite difference version of the right hand side of equation (\ref{mcf2}) with the spatial derivatives evaluated using equations (\ref{fd3}) and (\ref{fd4}).  The standard von Neumann stability analysis of equation (\ref{Euler2}) reveals that the time step must satisfy the Courant condition
\begin{equation}
\Delta t < {\textstyle {\frac 1 2}} {{(\Delta r)}^2}
\label{CFL2}
\end{equation}
The approximations in equations (\ref{fd3}) and (\ref{fd4}) are good at the tip.  However, far from the tip the surface asymptotically approaches a cylinder, which means that ${z_r} \to \infty$; consequently, equations \eqref{fd3} and \eqref{fd4} fail to be accurate in this region.  Thus, we have one equation (\ref{mcf1}) which does not work near the tip, and another (\ref{mcf2}) which does not work in the asymptotic region.

Our solution is to use both equations, with equation (\ref{mcf1}) used in a region that does not include the tip, and equation (\ref{mcf2}) used in a region that does not include the asymptotic region.  Often in numerical methods with two different grids, one makes the grids meet at a fixed boundary point.  Instead, in keeping with the spirit of differential geometry, we make our grids overlap, like two coordinate patches in an atlas.  In a numerical simulation, the last point in a grid cannot be evolved using the same method as the interior points and must be specified in some other way, typically some sort of boundary condition.  However, in our overlap method the final point of one grid corresponds to an interior point of the other grid. Thus, we determine the evolution of the finite difference function values at the end points of each grid using the function values at the corresponding interior points of the alternate grid.

Specifically, consider the grid values $z_i$ representing the function $z(r)$.  Let $N$ be the maximum value of $i$; we must specify the value of $z_N$ corresponding to the value that $z$ takes if $r=(N-1)\Delta r$ (since the grid points have spacing $\Delta r$ and the first grid point is at $r=0$).  We also have the grid values $r_i$ corresponding to the function $r(z)$ represented by the values on the other grid. If $r=(N-1)\Delta r$ is contained between some $r_i$ and $r_{i+1}$, then we interpolate linearly between the values of $z$ on that grid to obtain an interpolated value $z_N$.  If $r=(N-1)\Delta r$ is not contained between any $r_i$ on the other grid (as is sometimes the case as the surface moves under MCF), then we simply remove this last point from the grid, and subsequently the next to last grid point becomes the last one. Correspondingly, this method is used to specify the value of $r$ at the endpoint of the other grid, and to remove that endpoint if necessary.  Endpoints can also be removed from the grid if derivatives become so large at those points that they cannot be accurately computed using finite differences.

We end this section by explaining how we choose the perturbed initial data sets for the numerical simulations. We first recall that for an unperturbed solution constructed in \cite{IW19} (or \cite{IWZ20}, resp.), its initial data set  is obtained by joining a scaled bowl soliton to a cylinder at spatial infinity, and is defined precisely in the rescaled coordinates by equation \eqref{id} (or \eqref{id-crit}, resp.). Expressed in terms of the unscaled $(z,r)$-coordinates, the ODE for the bowl soliton takes the form
\begin{equation}\label{bowl-ode}
z_{rr} = (1+z_r^2)\left(\beta - \frac{z_r}{r}\right).
\end{equation}

For each choice of $\gamma$, we fix $c=1$ and choose $\tau_0$ large (in all the numerical simulations we let $\tau_0=4$) so that the matched asymptotics explained in Section \ref{setup} make sense. We then choose $\beta$ according to $\beta=c^{-1}(\gamma-1/2)2^{-(\gamma-1/2)} e^{-(\gamma+1/2)\tau_0}$, and numerically integrate equation \eqref{bowl-ode} outward from $r=0$ to value $r_1$, which is obtained by writing $R_1$ (chosen to be $e^{\gamma \tau_0/2}$ in our numerical simulations) in the unscaled $r$-coordinate. For all $r>r_1$, we use the following analytic formula which follows from rewriting the second equation in \eqref{id} in the $(r,z)$-coordinates:
\begin{equation}\label{formula}
z(r) = z(r_1)+\frac{1}{c}\left[ \left( 2e^{-\tau_0}-r^2\right) ^{\frac{1}{2}-\gamma} - \left(2e^{-\tau_0}-r_1^2\right)^{\frac{1}{2}-\gamma} \right].
\end{equation}
In the patch where we write $z$ as a function of $r$, we use formula \eqref{formula}. In the patch where $r$ is written as a function of $z$, we use the formula obtained by inverting the expression \eqref{formula}.

To implement the perturbations of the initial data sets that we consider in this work, we distort the above initial data set by putting a ``dimple'' near the tip (for the near class, cf. Figures \ref{near1}--\ref{near4} in Section \ref{near_class}) or putting a dimple in the region far away from the tip (for the far class, cf. Figures \ref{far1}--\ref{far3} in Section \ref{far_class}). Precisely, for a perturbed initial data set expected to be in the near class, we choose an amplitude $a_0$ and a maximum $r$ value $r_m$. We leave the surface undistorted for $r>r_m$, while for $r\leq r_m$ we change $z$ as follows,
\begin{align*}
z\mapsto z + a_0\left(1-\frac{r^2}{r_m^2} \right)^2.
\end{align*}
For a perturbed initial data set expected to be in the far class, we choose an amplitude $a_0$, and two values of $z$: $z_a$ and $z_b$. We distort the surface only on the interval $(z_a, z_b)$ by letting
\begin{align*}
r\mapsto r-a_0\frac{(z-a_a)^2(z-z_b)^2}{(z_b-z_a)^4}.
\end{align*}

\section{Numerical results}\label{num_results}

Let $\Gamma_t$, where $t\in[t_0,T)$, be a MCF solution of a noncompact hypersurface constructed in \cite{IW19, IWZ20}. In particular, $\Gamma_{t}$ is asymptotic to a shrinking cylinder $\cyl_t$ near spatial infinity for all $t\in[t_0, T)$, where $T$ is the time precisely when $\cyl_t$ collapses to a line \cite{SS14} (see also Figure \ref{fig1}). At the tip of $\Gamma_t$, the curvature blows up at a Type-II rate $(T-t)^{-(\gamma+1/2)}$ for $\gamma\geq 1/2$; at spatial infinity, $\Gamma_t$ is asymptotic to $\cyl_t$, which is a self-shrinking solution to MCF. We regard $\Gamma_t$ as MCF with a degenerate neckpinch forming at spatial infinity.

Let $\tilde\Gamma_{t_0}$ be a perturbation of $\Gamma_{t_0}$ such that $\tilde\Gamma_{t_0}$ is still an embedding of a noncompact hypersurface and asymptotic to $\cyl_{t_0}$. Let $r(z)$ be the profile of rotation for $\Gamma_{t_0}$ and $\tilde r(z)$ the profile of rotation for $\tilde{\Gamma}_{t_0}$. In particular, by construction in \cite{IW19, IWZ20}, $r'(z)$ is always positive.

We now present the results from numerical simulations by solving MCF with initial data $\tilde\Gamma_{t_0}$ using the overlapping method explained in Section \ref{num_method}.

\subsection{The near class}\label{near_class}

We define an initial data set $\tilde\Gamma_0$ to be contained in the \emph{near} class if the MCF solution starting from $\tilde\Gamma_0$ develops a Type-II singularity (modelled by the bowl soliton) at $t=\tilde T\leq T$. Our numerical simulations show that if the perturbation of $r(z)$ is placed sufficiently close to the tip so that $\tilde\Gamma_{t_0}$ is still locally a graph over the cross-sectional ball of $\cyl_{t_0}$, then the evolution of $\tilde\Gamma_t$ starting from $\tilde\Gamma_{t_0}$ generated by $\tilde{r}(z)$ resembles that of $\Gamma_t$.

\begin{figure}[H]
	\centering
	\begin{subfigure}[H]{0.8\linewidth}
		\includegraphics[width=\linewidth]{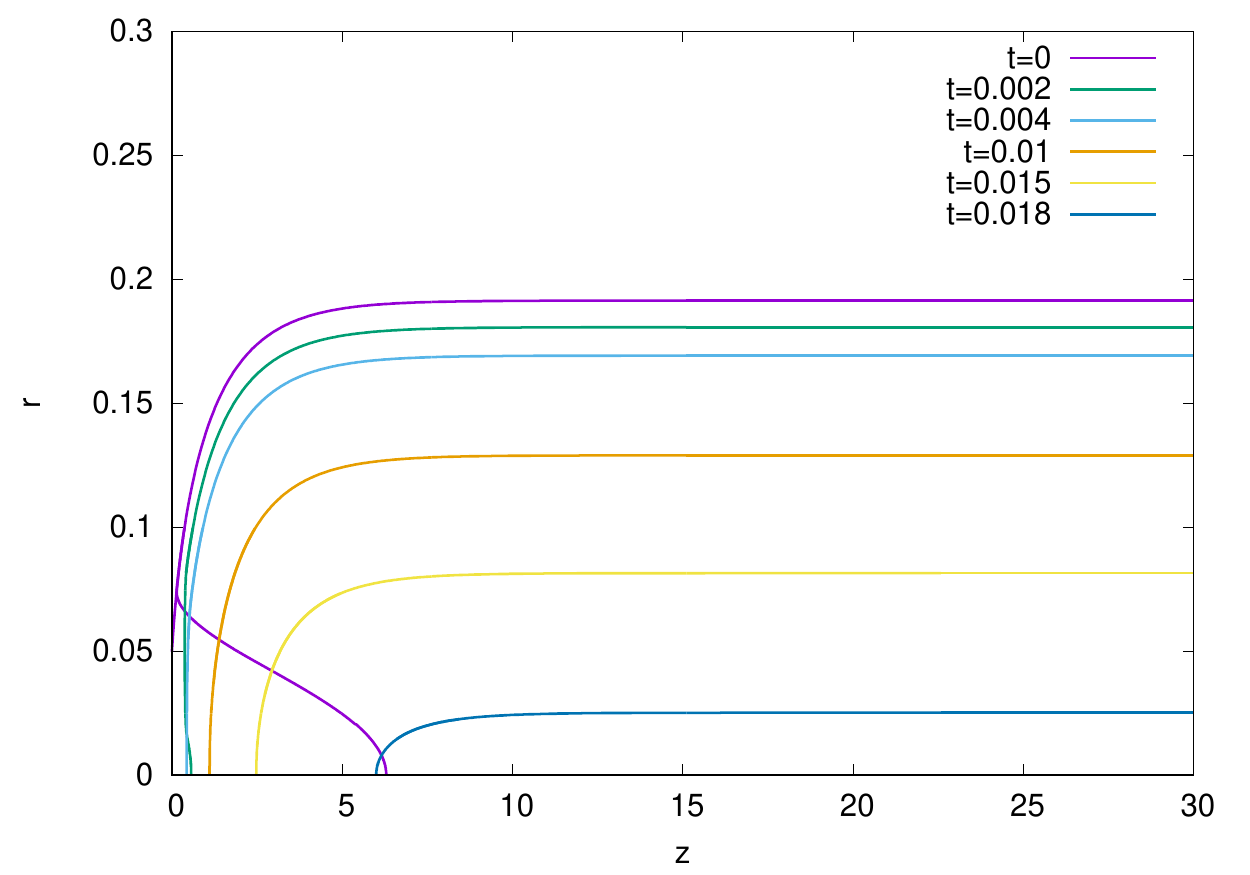}
		\caption{MCF evolutions}\label{nearshapeg-p5}
	\end{subfigure}\vspace{12pt}
	\begin{subfigure}[H]{0.4\linewidth}
		\includegraphics[width=\linewidth]{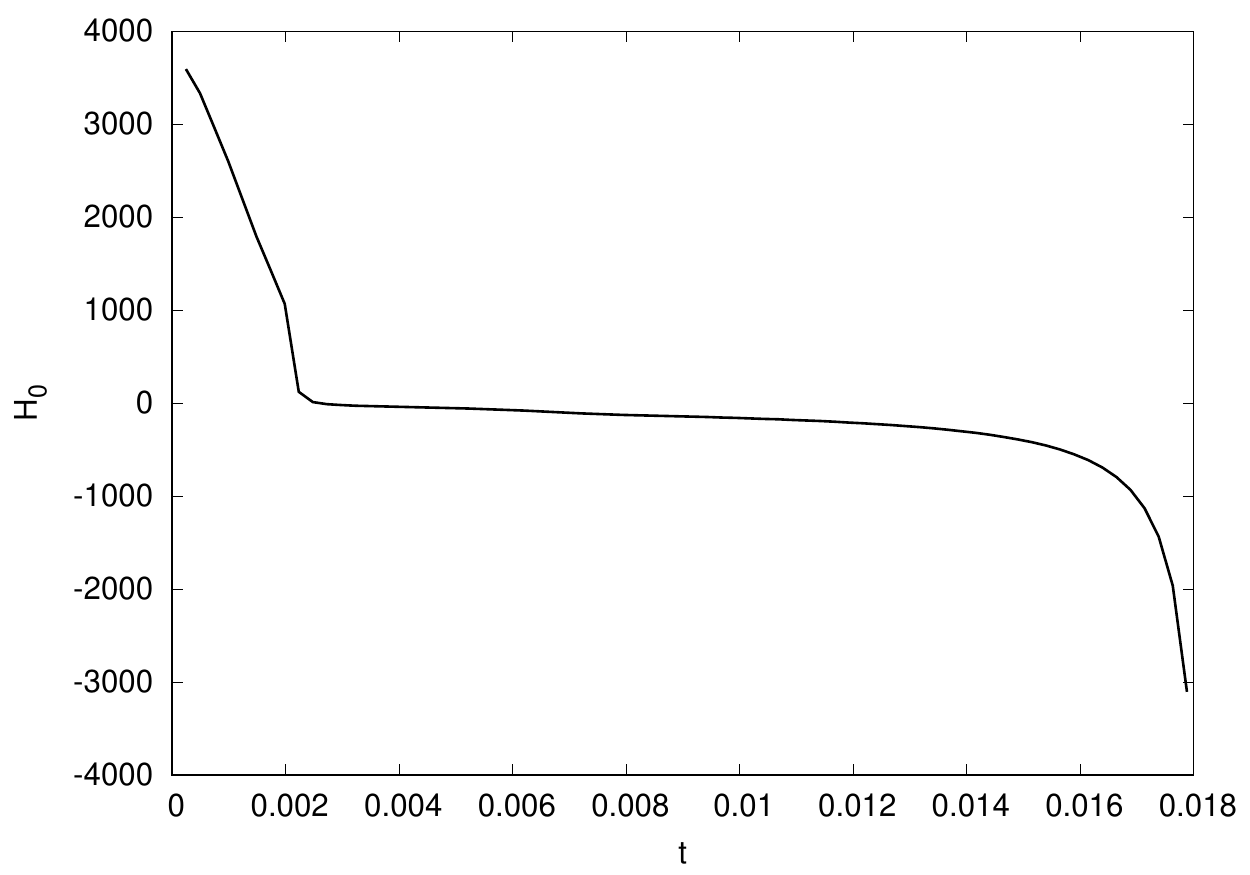}
		\caption{Unscaled curvature at the tip}\label{nearcurvg-p5}
	\end{subfigure}
	\begin{subfigure}[H]{0.4\linewidth}
		\includegraphics[width=\linewidth]{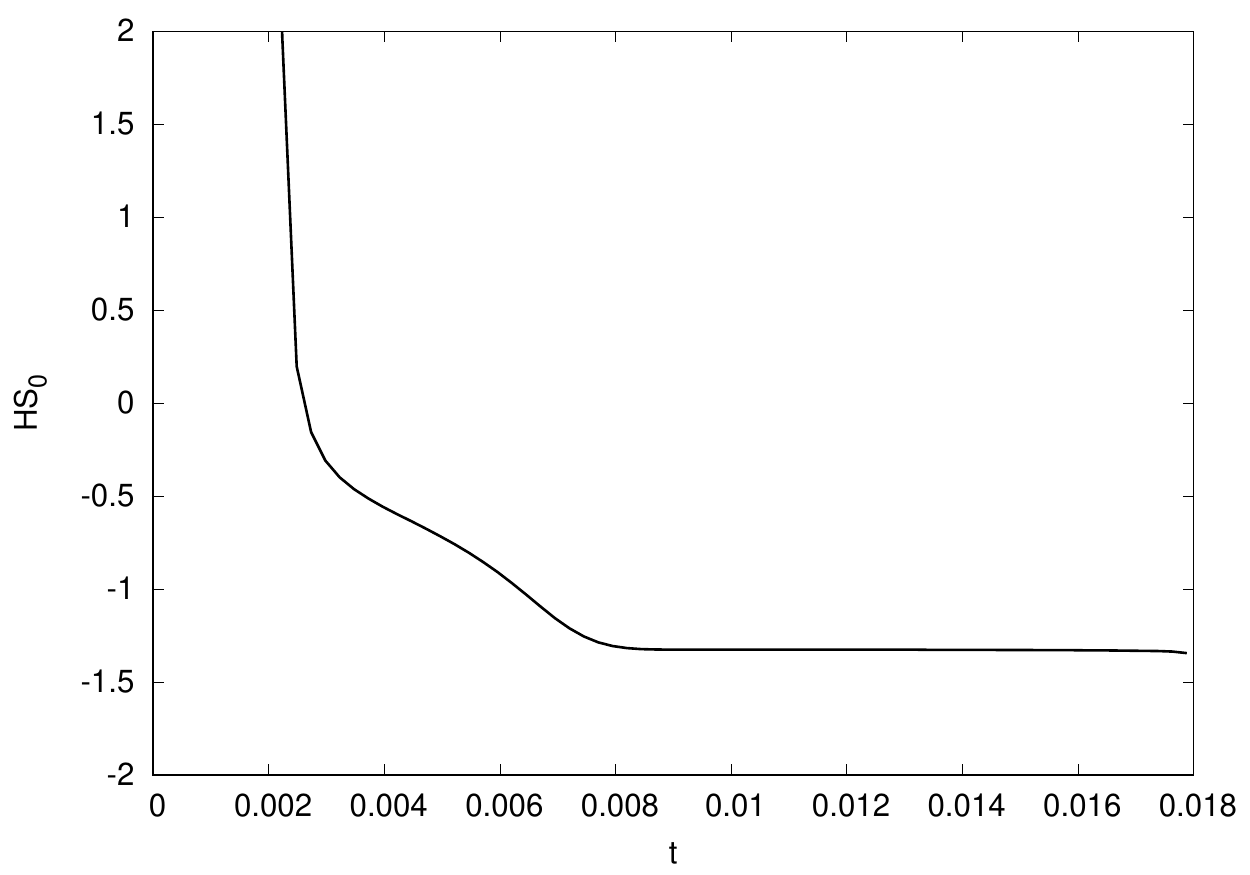}
		\caption{Rescaled curvature at the tip}\label{nearcurv2g-p5}
	\end{subfigure}
	\caption{Numerical simulation of a MCF solution in the near class, for $\gamma=1/2$.}\label{near1}
\end{figure}

\begin{figure}[H]
	\centering
	\begin{subfigure}[H]{0.8\linewidth}
		\includegraphics[width=\linewidth]{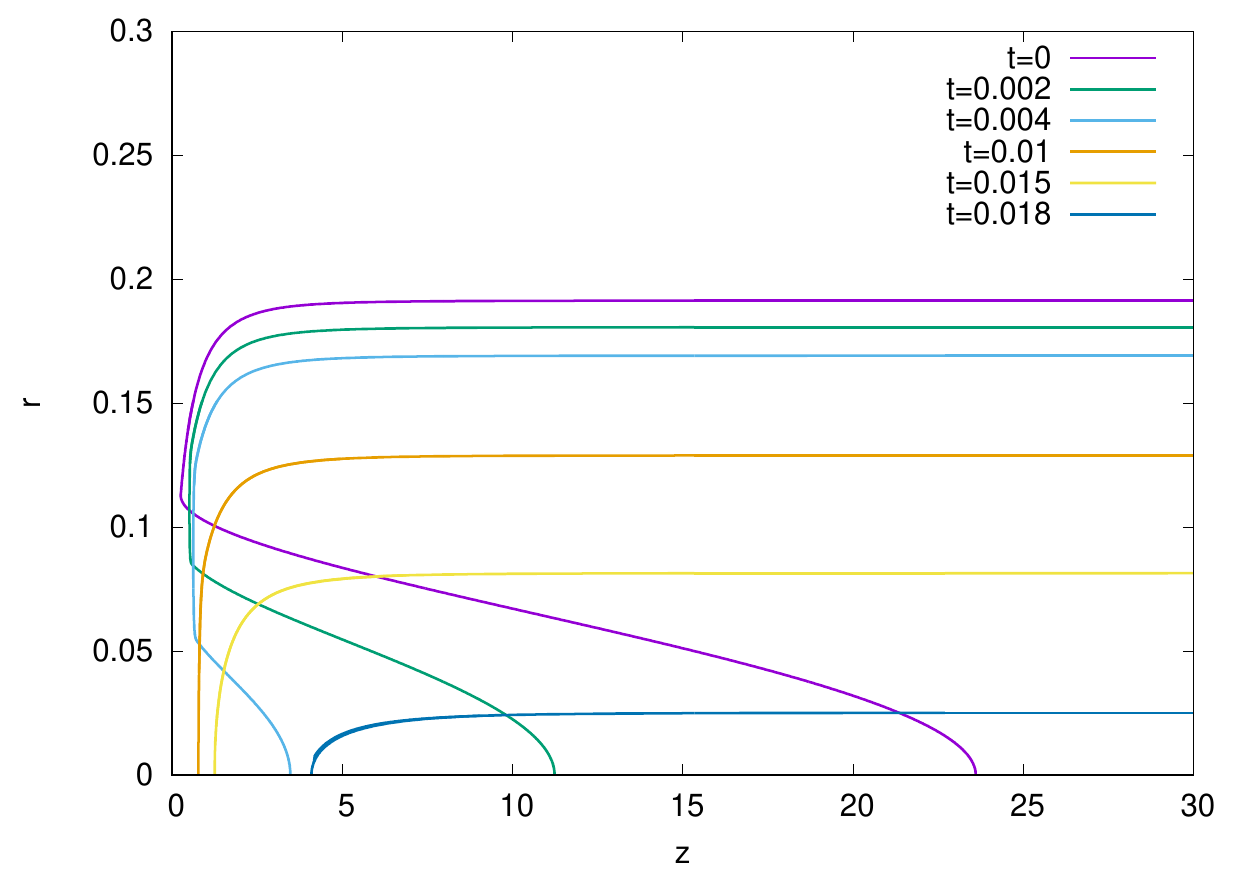}
		\caption{MCF evolutions}\label{nearshapeg-p75}
	\end{subfigure}\vspace{12pt}
	\begin{subfigure}[H]{0.4\linewidth}
		\includegraphics[width=\linewidth]{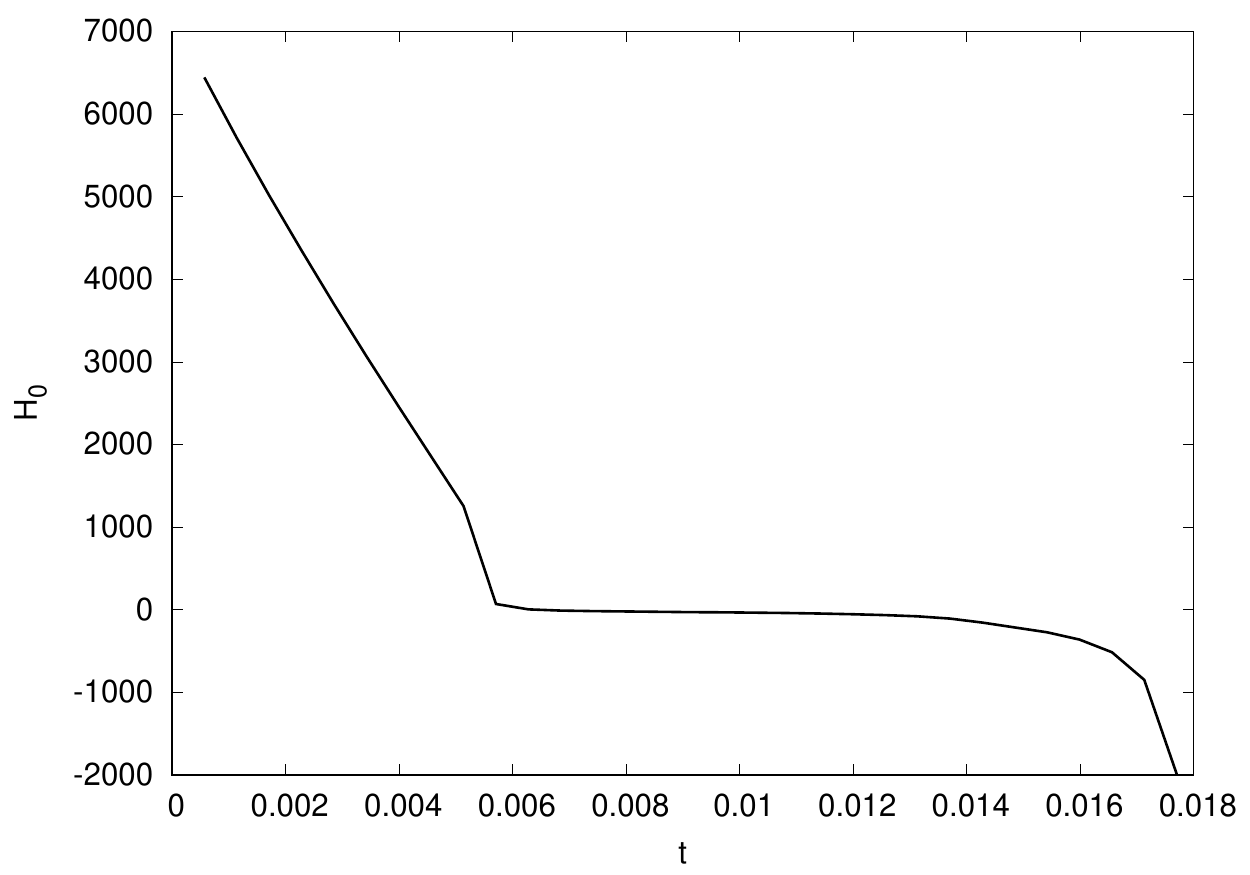}
		\caption{Unscaled curvature at the tip}\label{nearcurvg-p75}
	\end{subfigure}
	\begin{subfigure}[H]{0.4\linewidth}
		\includegraphics[width=\linewidth]{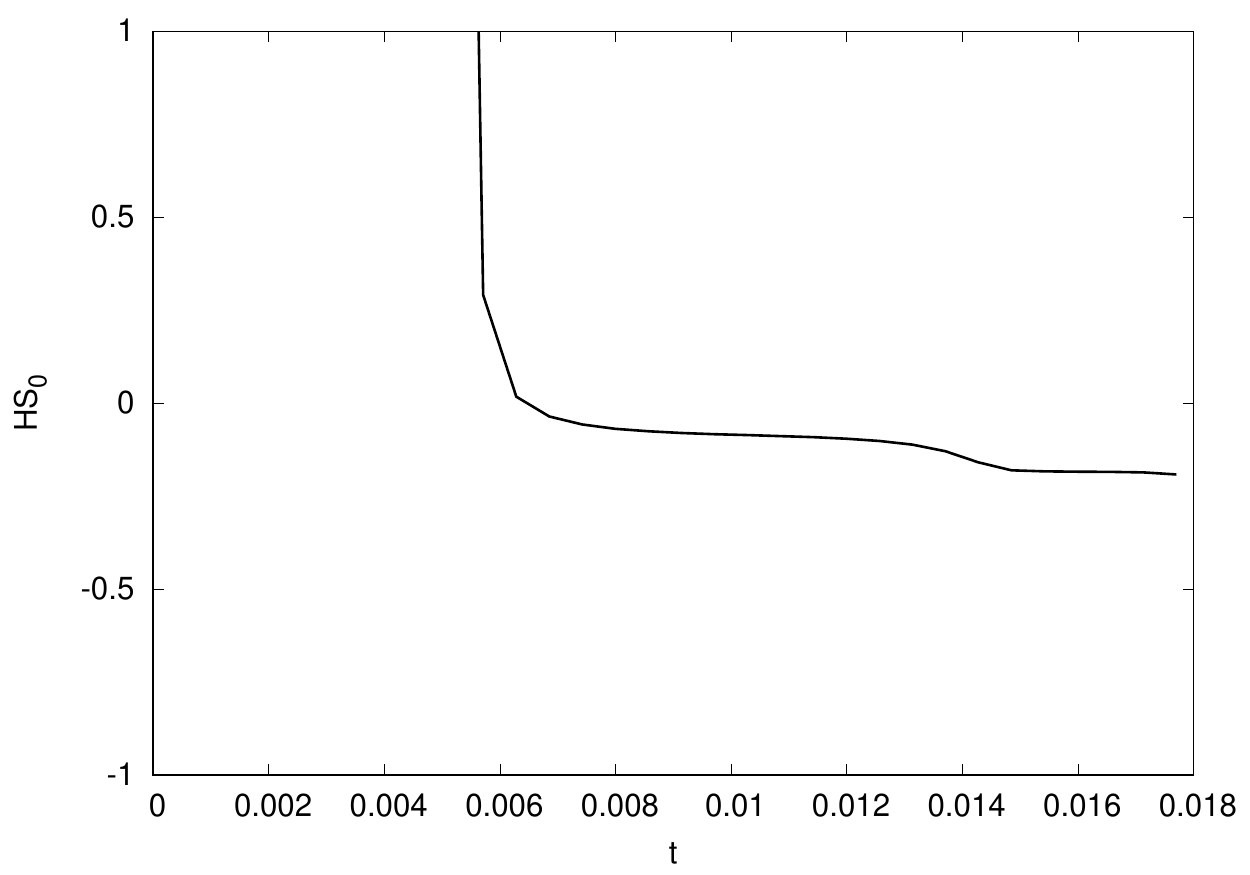}
		\caption{Rescaled curvature at the tip}\label{nearcurv2g-p75}
	\end{subfigure}
	\caption{Numerical simulation of a MCF solution in the near class, for $\gamma=3/4$.}\label{near2}
\end{figure}

\begin{figure}[H]
	\centering
	\begin{subfigure}[H]{0.8\linewidth}
		\includegraphics[width=\linewidth]{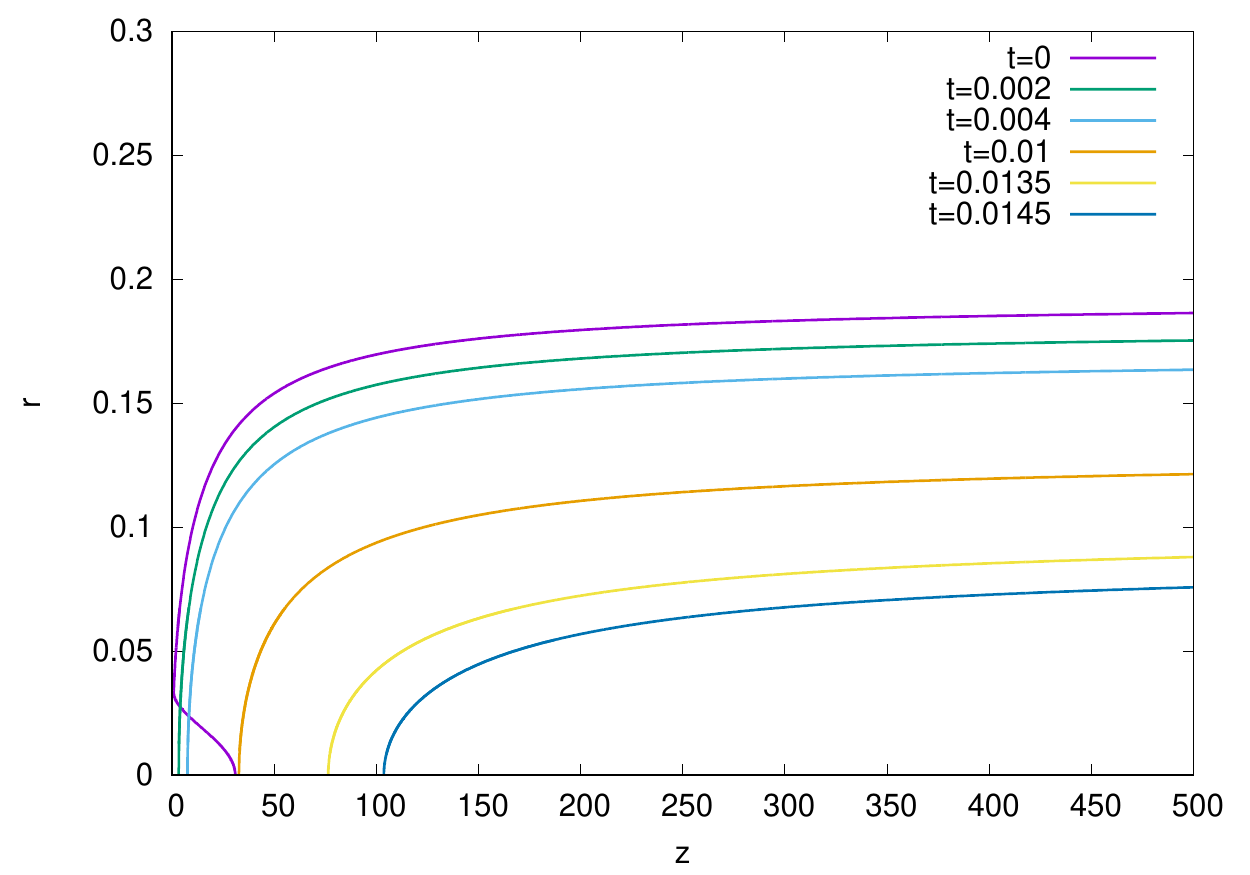}
		\caption{MCF evolutions}\label{nearshapeg-1p5}
	\end{subfigure}\vspace{12pt}
	\begin{subfigure}[H]{0.4\linewidth}
		\includegraphics[width=\linewidth]{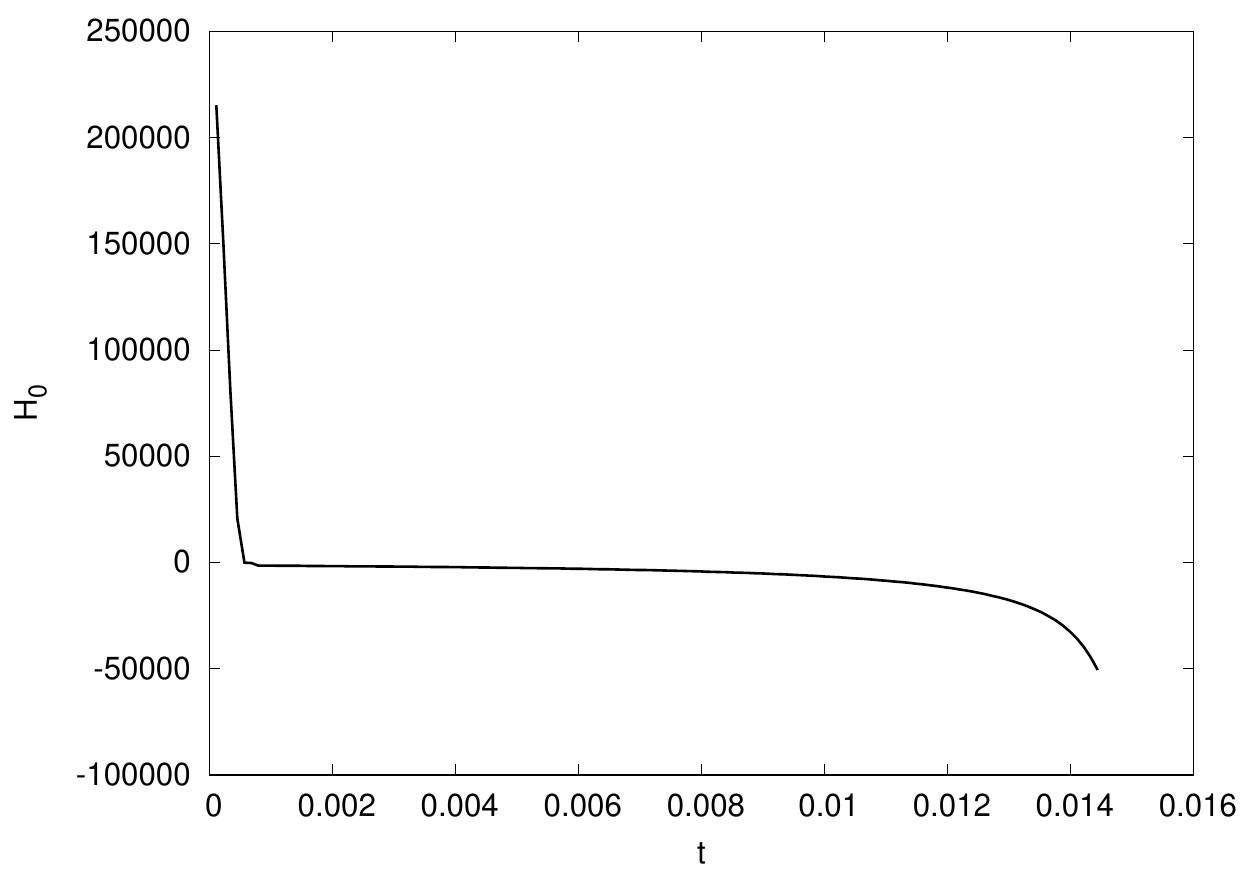}
		\caption{Unscaled curvature at the tip}\label{nearcurvg-1p5}
	\end{subfigure}
	\begin{subfigure}[H]{0.4\linewidth}
		\includegraphics[width=\linewidth]{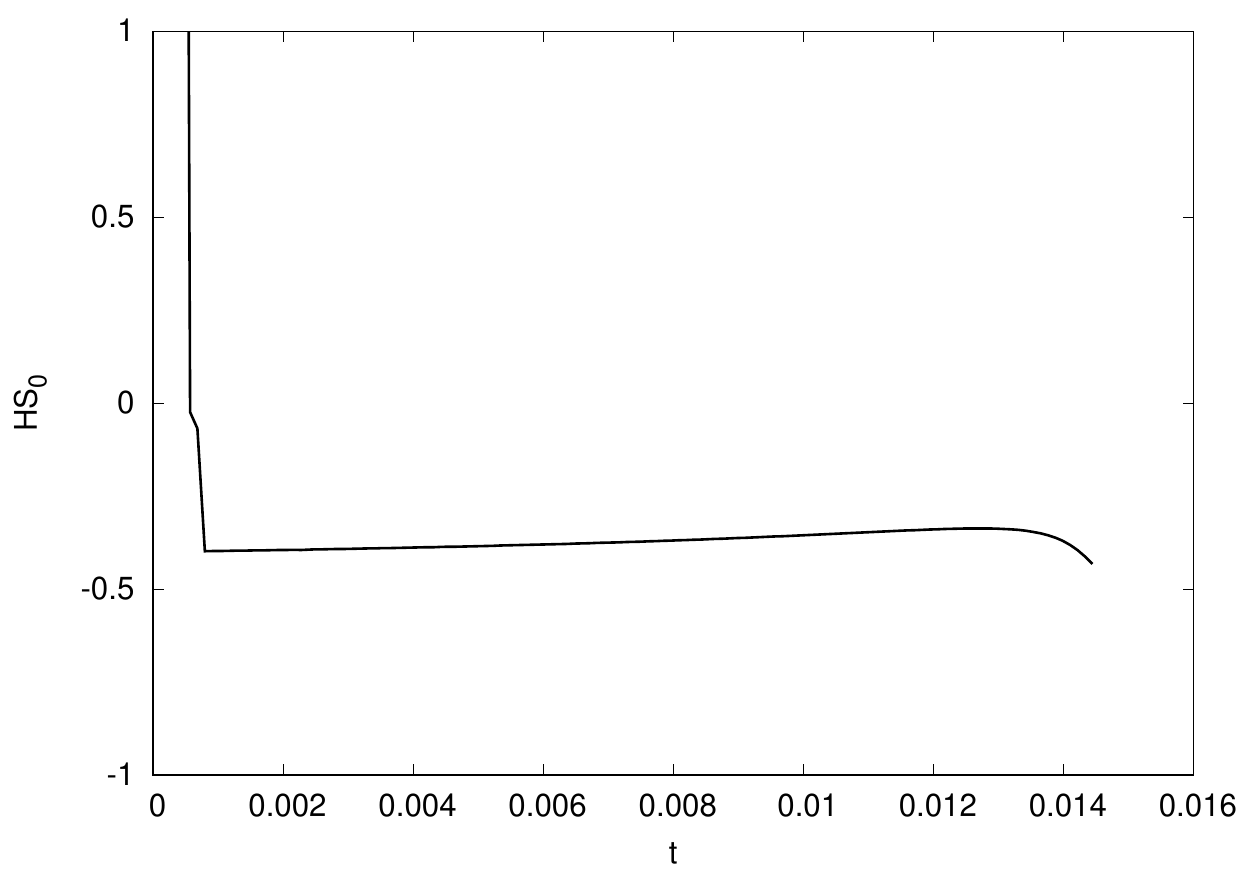}
		\caption{Rescaled curvature at the tip}\label{nearcurv2g-1p5}
	\end{subfigure}
	\caption{Numerical simulation of a MCF solution in the near class, with the initial perturbation away from the tip, for $\gamma=3/2$.}\label{near3}
\end{figure}

\begin{figure}[H]
	\centering
	\begin{subfigure}[H]{0.8\linewidth}
		\includegraphics[width=\linewidth]{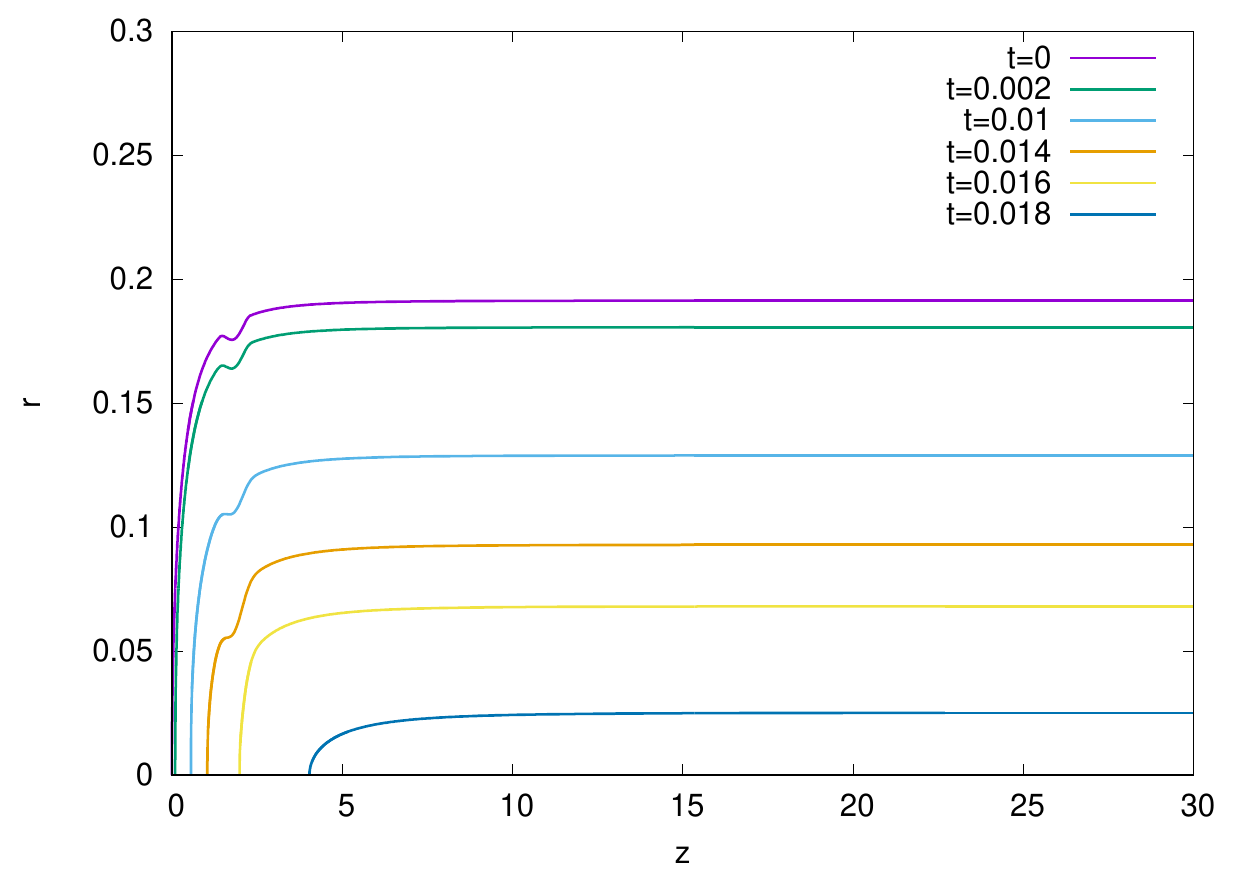}
		\caption{MCF evolutions}\label{subcritshape}
	\end{subfigure}\vspace{12pt}
	\begin{subfigure}[H]{0.4\linewidth}
		\includegraphics[width=\linewidth]{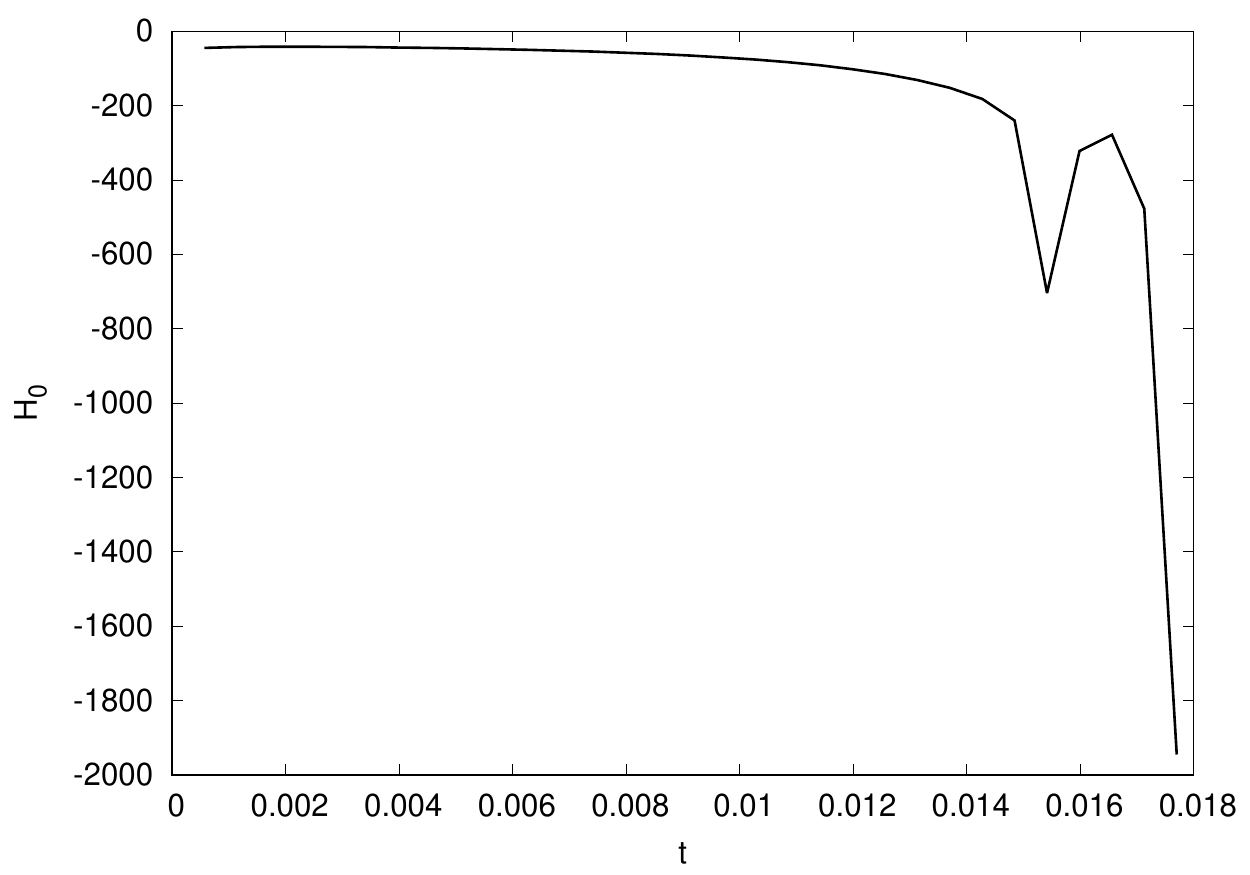}
		\caption{Unscaled curvature at the tip}\label{subcritcurv}
	\end{subfigure}
	\begin{subfigure}[H]{0.4\linewidth}
		\includegraphics[width=\linewidth]{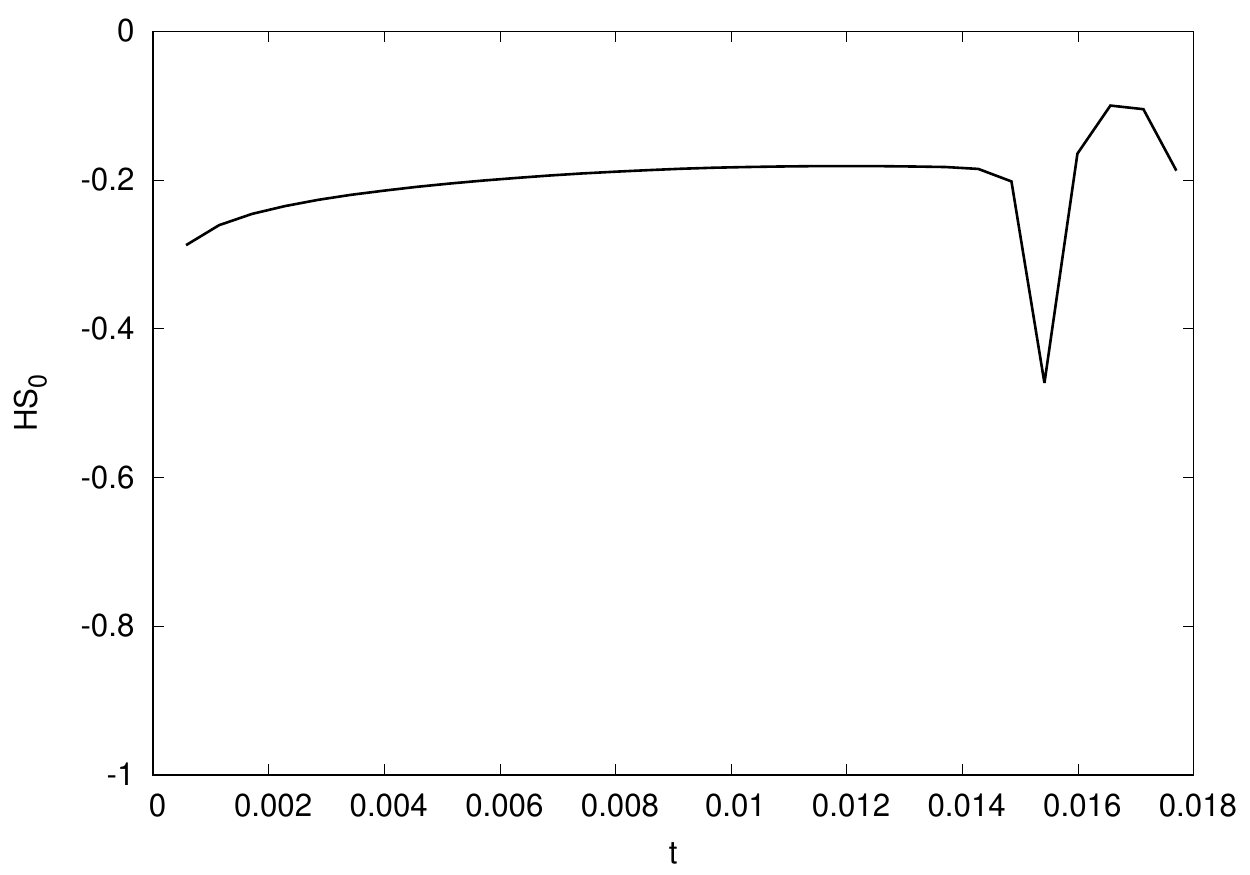}
		\caption{Rescaled curvature at the tip}\label{subcritcurv2}
	\end{subfigure}
	\caption{Numerical simulation of a MCF solution in the near class, for $\gamma=3/4$.}\label{near4}
\end{figure}

We now present some numerical simulations of MCF originating in the near classes. See Figure \ref{near1} for a simulation with $\gamma=1/2$, and Figures \ref{near2} and \ref{near3} for simulations with two different values of $\gamma>1/2$. The perturbation of the initial data set in each of Figures \ref{near1}--\ref{near3} involves a dimple placed very close to the tip. The near class also includes initial data sets with the dimple placed relatively far from the tip, as shown in Figure \ref{near4}. In each figure, the unscaled (signed) curvature at the tip is denoted by $H_0$, whereas $HS_0:= H_0(T-t)^{\gamma+1/2}$ denotes the rescaled curvature at the tip according to the Type-II rate of the unperturbed solution $\Gamma_t$.

Figure \ref{near1} shows the numerical simulation of a MCF solution in the near class, for $\gamma=1/2$. Graph (A) shows the embeddings as time progresses. Note that the perturbation in the initial data near the tip quickly disappears. Graph (B) shows the time development of the mean curvature at the tip. Note that it becomes arbitrarily large near the singularity. Graph (C) shows the time development of the rescaled mean curvature at the tip. Note that it approaches a constant, consistent with a Type-II singularity forming at the tip.

Figure \ref{near2} shows the numerical simulation of a MCF solution in the near class, for $\gamma=3/4$.
Graph (A) again shows the time progression of the embeddings; in this case the perturbation in the initial data disappears a bit more slowly. Graph (B) shows the time development of the mean curvature at the tip again becoming arbitrarily large, while Graph (C) shows the rescaled mean curvature at the tip approaching a constant, consistent with a Type-II singularity.

Figure \ref{near3} shows the numerical simulation of a MCF solution in the near class, for $\gamma=3/2$. Graphs (A)--(C) indicate behavior very similar to that seen in the examples illustrated in Figure \ref{near1} and Figure \ref{near2}.

Figure \ref{near4} shows the numerical simulation of a MCF solution in the near class with the perturbation of the initial data set located away from the tip, for $\gamma=3/4$. The bahaviors indicated in Graphs (A)--(C) resemble those observed in the examples illustrated in Figures \ref{near1}--\ref{near3}.

According to these numerical simulations, the mean curvature flow $\tilde\Gamma_t$ developing from initial data in the near class develops a Type-II singularity at the tip for a wide range of choices of $\gamma$. Moreover, the above figures of the rescaled curvature at the tip suggest that the Type-II singularity of $\tilde\Gamma_t$ blows up at the same Type-II rate. Thus, we are led to the following conjecture: 
\begin{conj}
	Mean curvature flow solutions starting from initial data sets in the near class have the same singular behavior as the unperturbed solution, their maximum curvatures blow up at the same Type-II rate, and they always converge to a corresponding unperturbed solution.
\end{conj}

\begin{remark}
	In \cite{CSS07}, Clutterbuck, Schn\"urer and Schulze prove the stability of the bowl soliton by showing that the MCF of an entire graph that is asymptotic to the bowl soliton converges uniformly to the bowl soliton as $t\to\infty$. In comparison, since $\Gamma_{t_0}$ is close to the bowl soliton in the interior region, and so is the perturbed initial data $\tilde{\Gamma}_{t_0}$, our numerical results suggest a localized stability result for the bowl soliton.
\end{remark}

\begin{remark}
	We recall the ``vertical line test'' which is a consequence of the Sturmian theorem \cite[Section 4]{AAG} and says that the number of intersections of the graph of $r(z,t)$ and a vertical line $z=z_0$ is nonincreasing in time. Our numerical simulation shows that a solution $\tilde\Gamma_t$ in the near class converges to a solution $\Gamma_t$ from \cite{IW19, IWZ20}. In particular, this is consistent with the vertical line test as any intersection at the initial time disappears as the graph shrinks and moves toward spatial infinity under mean curvature flow. The vertical line test does not imply any convergence of the solution. In comparison, our numerical simulation indicates the convergence of $\tilde\Gamma_t$ and the local stability of $\Gamma_t$ from \cite{IW19, IWZ20}.
\end{remark}

\subsection{The far class}\label{far_class}

We define an initial data set $\tilde\Gamma_{t_0}$ to be contained in the \emph{far} class if the MCF solution starting from $\tilde\Gamma_{t_0}$ develops a local Type-I neckpinch singularity (which is modelled by the cylinder $S^{n-1}\times\mathbb{R}$) at $t=\tilde T\leq T$. Our numerical simulations show that if the perturbed initial profile $\tilde{r}(z)$ contains a local minimum --- no matter the size of its depth --- sufficiently far from the tip region of $\Gamma_{t_0}$ (equivalently, in the region where $|r'(z)|$ is sufficiently small), then MCF starting from $\tilde{\Gamma}_{t_0}$ generated by $\tilde r(z)$ develops a Type-I nondegenerate neckpinch.

In the literature, Type-I singularities are often called \emph{rapidly forming,} because one has a bound on the  time remaining until the singularity form of the type $\sqrt{T-t}\leq C/(\sup_{\Gamma_t}|h|)$; by contrast, Type-II singularities are called \emph{slowly forming} because no such bound holds. Our results for perturbations in the far class illustrate this heuristically: a Type-II singularity is forming slowly at the tip. However, the solution does not have time to become singular there, because a Type-I singularity is rapidly forming away from the tip.

\begin{figure}[H]
	\centering
	\includegraphics[width=0.65\linewidth]{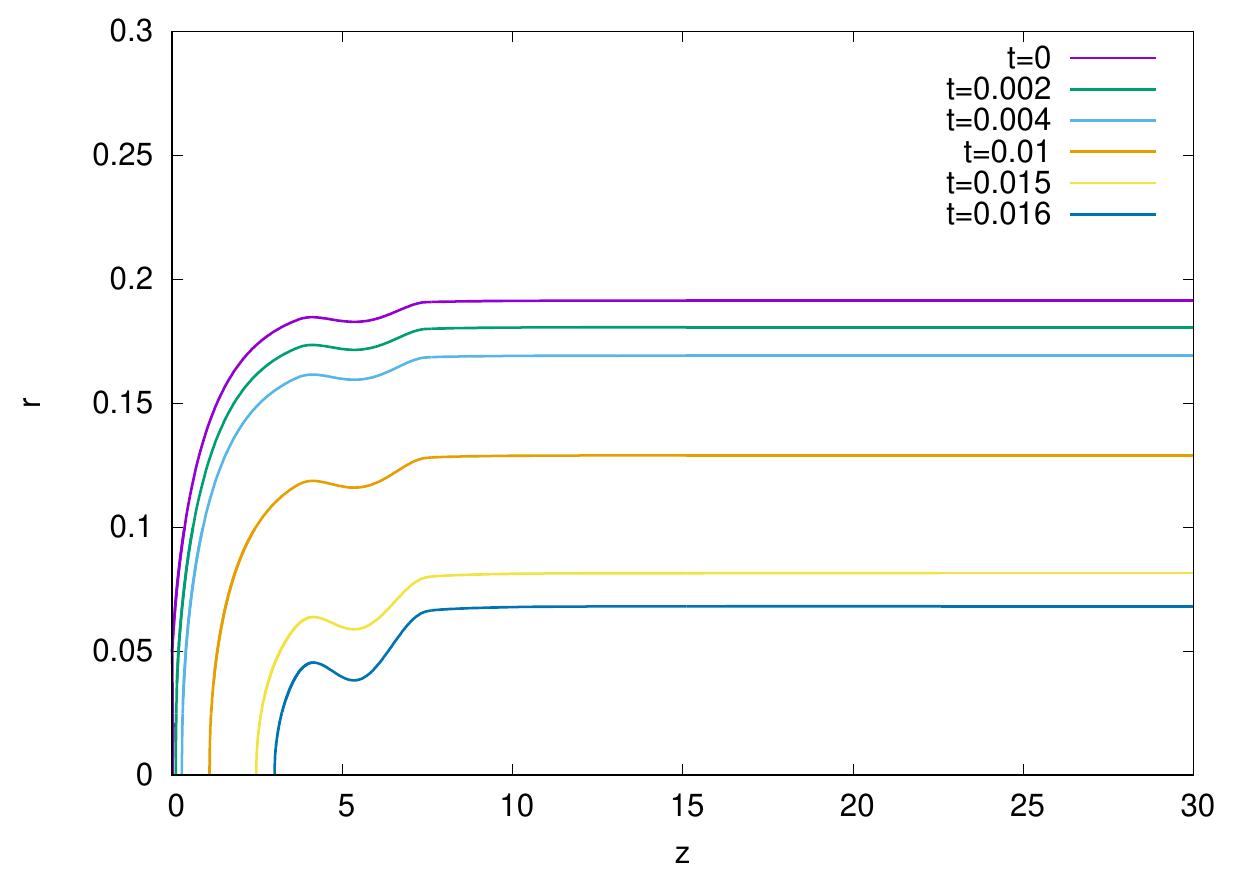}
	\caption{Numerical simulation of a MCF solution in the far class, for $\gamma=1/2$.}\label{far1}
\end{figure}

\begin{figure}[H]
	\centering
	\includegraphics[width=0.65\linewidth]{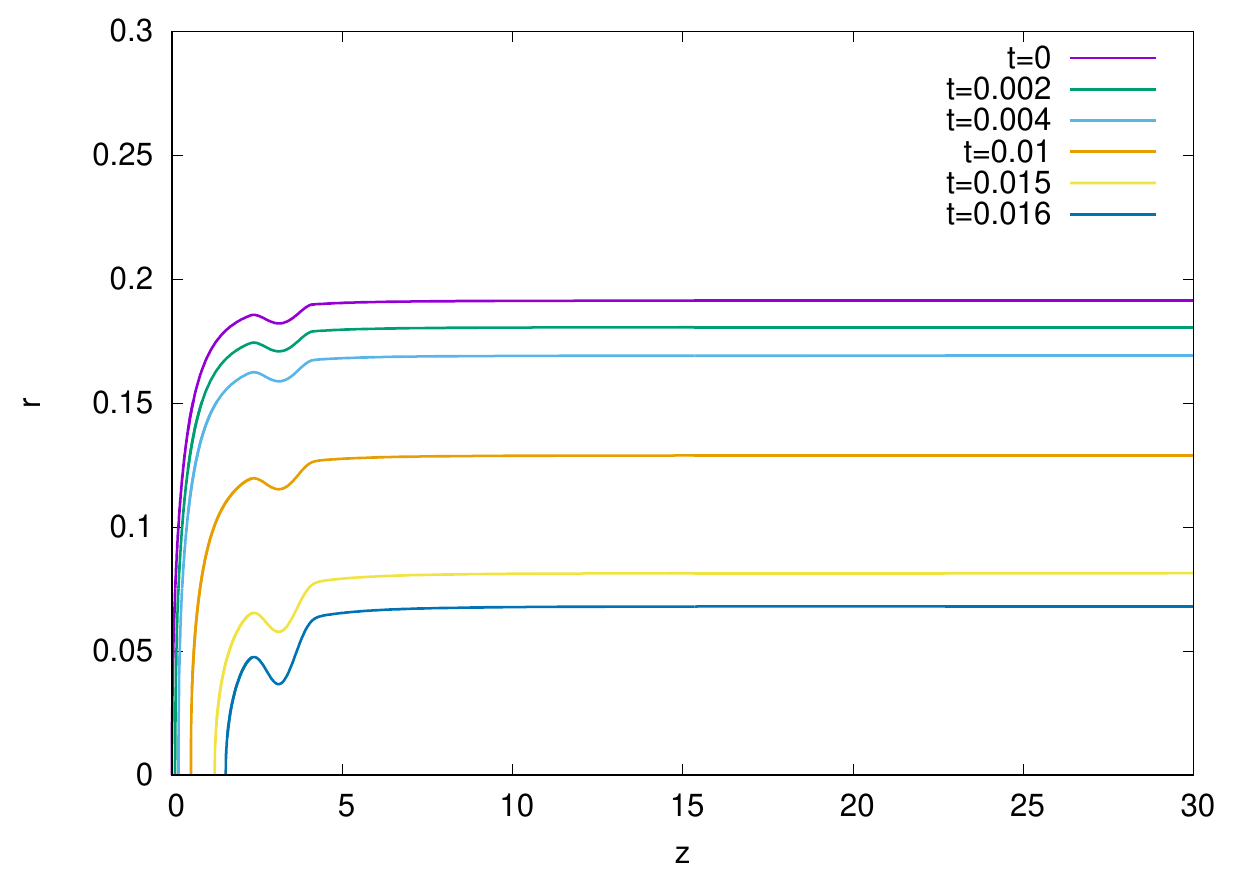}
	\caption{Numerical simulation of a MCF solution in the far class, for $\gamma=3/4$.}\label{far2}
\end{figure}

\begin{figure}[H]
	\centering
	\includegraphics[width=0.65\linewidth]{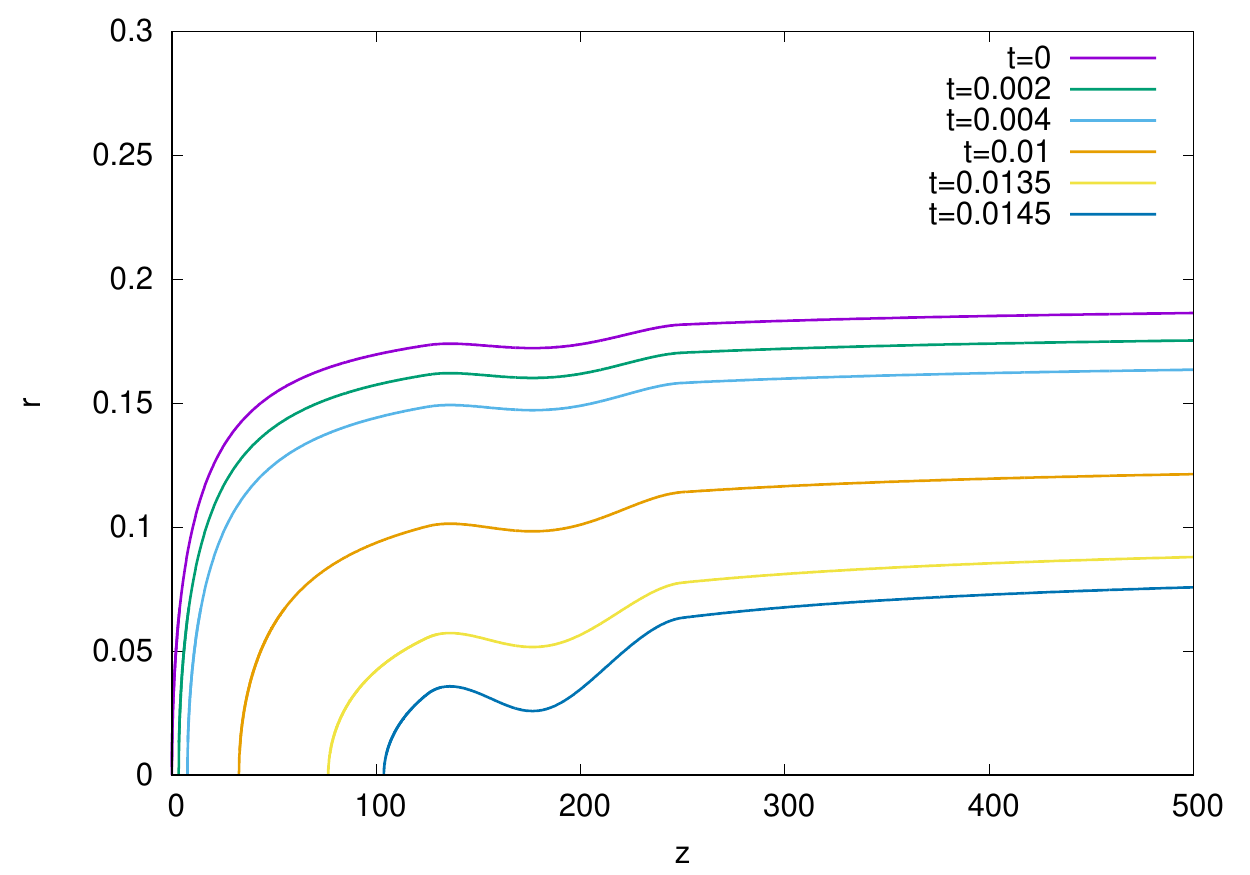}
	\caption{Numerical simulation of a MCF solution in the far class, for $\gamma=3/2$.}\label{far3}
\end{figure}

Below we display a few numerical simulations of MCF originating in the far classes. See Figure \ref{far1} for a simulation with $\gamma=1/2$, and Figures \ref{far2} and \ref{far3} for simulations with two different values of $\gamma>1/2$. In each figure, the unscaled curvature at the tip is denoted by $H_0$, and $HS_0:= H_0(T-t)^{\gamma+1/2}$ denotes the rescaled curvature at the tip according to the Type-II rate of the unperturbed solution $\Gamma_t$.

Figure \ref{far1} shows the numerical simulation of a MCF solution in the far class, for $\gamma=1/2$. Note that as time progresses, the perturbation far from the tip increases with time, indicating that a nondegenerate neckpinch necessarily forms far away, albeit at a finite distance, from the tip. This singular behavior is different from that of a degenerate neckpinch forming at spatial infinity in the unperturbed solution.

Figure \ref{far2} shows the numerical simulation of a MCF solution in the far class, for $\gamma=3/4$. Again, we see that the perturbation far from the tip grows with time. The evolution indicates the likely formation of a nondegenerate neckpinch singularity far away, but still at a finite distance, from the tip. This shows a singular behavior distinct from that of the unperturbed solution.

Figure \ref{far3} shows the numerical simulation of a MCF solution in the far class, for $\gamma=3/2$. The perturbation far from the tip --- even though it is initially quite small --- increases with time and should result in a nondegenerate neckpinch singularity at a finite distance away from the tip, contrasting with the degenerate neckpinch at infinity in the unperturbed solution.

\begin{remark}
	Under MCF, a solution in the far class always has one intersection with a vertical line. This is consistent with the vertical line test.
\end{remark}

\subsection{Near and far classes in higher dimensions}\label{n=3}
As explained in Section \ref{num_method}, the numerical results in higher dimensions should be the same since the lower order term in the evolution equation just has a different constant coefficient. Our simulations agree with this expectation, as illustrated by the following sample figures for $n=3$.

\begin{figure}[H]
	\centering
	\begin{subfigure}[H]{0.8\linewidth}
		\includegraphics[width=\linewidth]{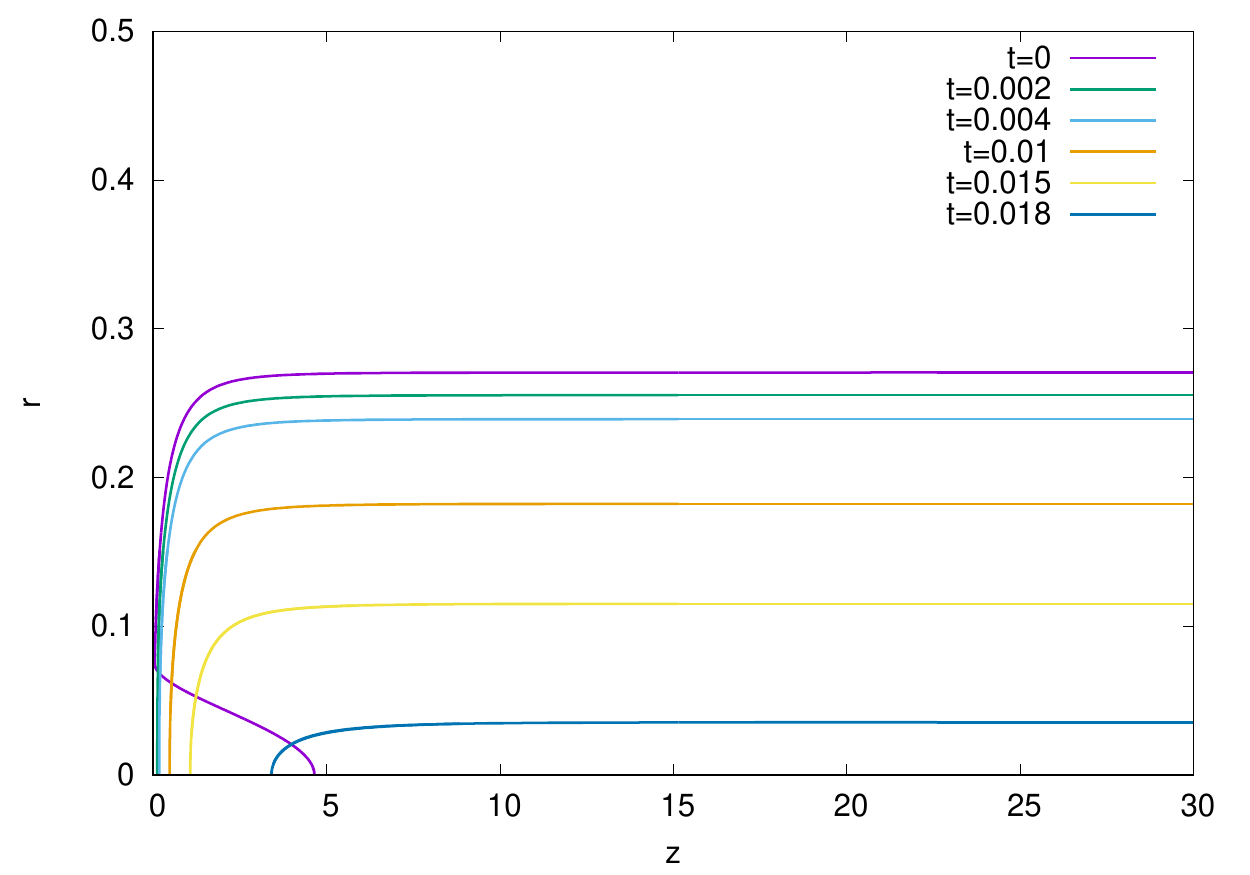}
		\caption{MCF evolutions}\label{nearshapegn3}
	\end{subfigure}\vspace{12pt}
	\begin{subfigure}[H]{0.4\linewidth}
		\includegraphics[width=\linewidth]{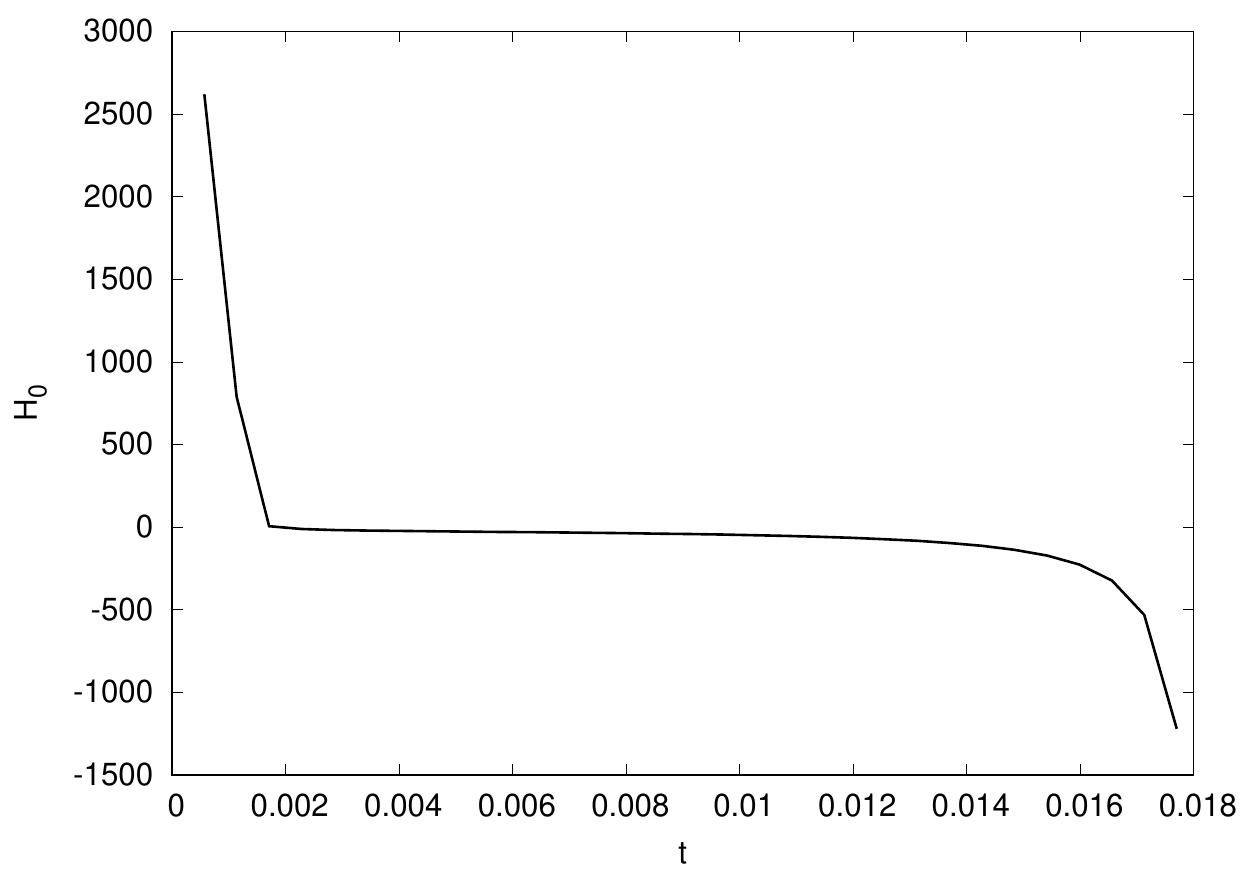}
		\caption{Unscaled curvature at the tip}\label{nearcurvn3}
	\end{subfigure}
	\begin{subfigure}[H]{0.4\linewidth}
		\includegraphics[width=\linewidth]{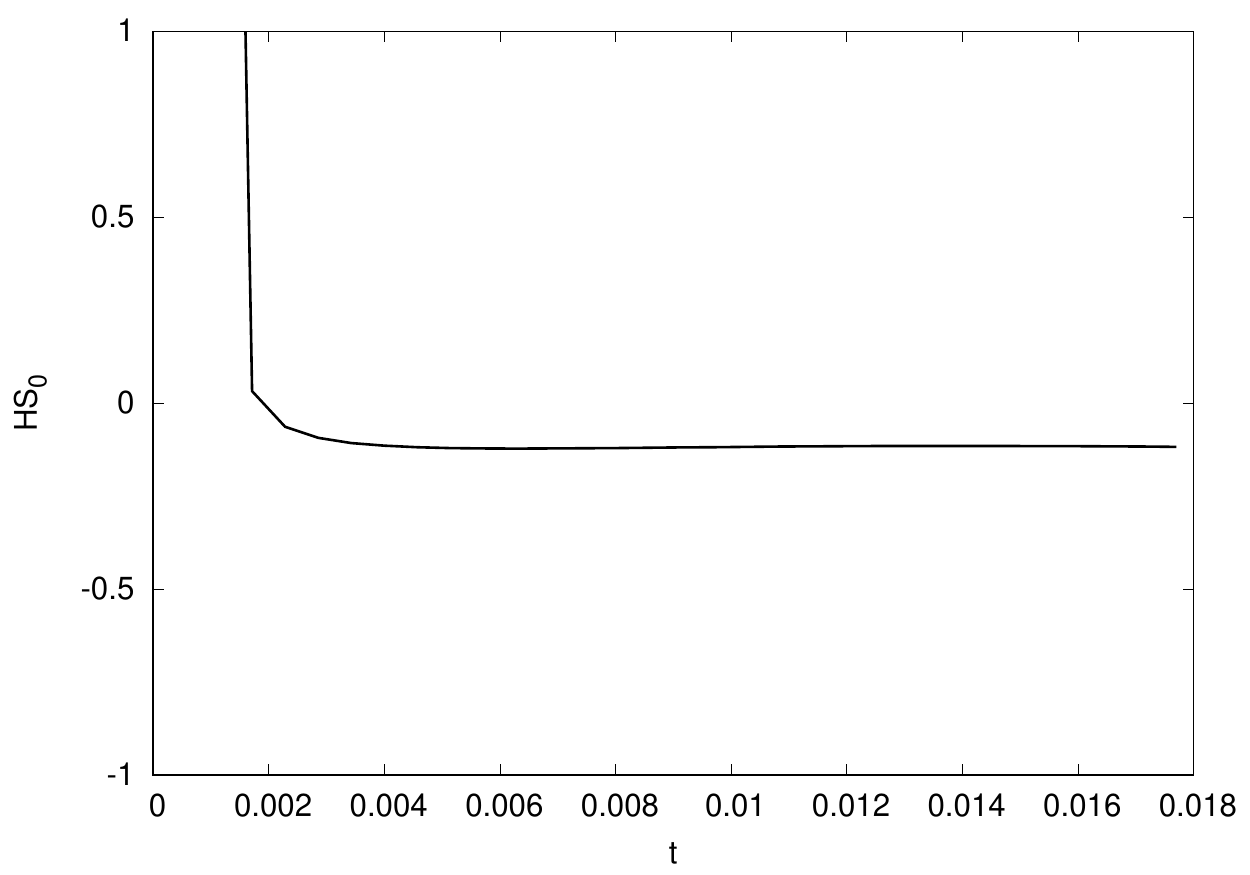}
		\caption{Rescaled curvature at the tip}\label{nearcurv2n3}
	\end{subfigure}
	\caption{Numerical simulation of a MCF solution in the near class, for $n=3$ and $\gamma=3/4$.}\label{near_n3}
\end{figure}

\begin{figure}[H]
	\centering
	\includegraphics[width=0.65\linewidth]{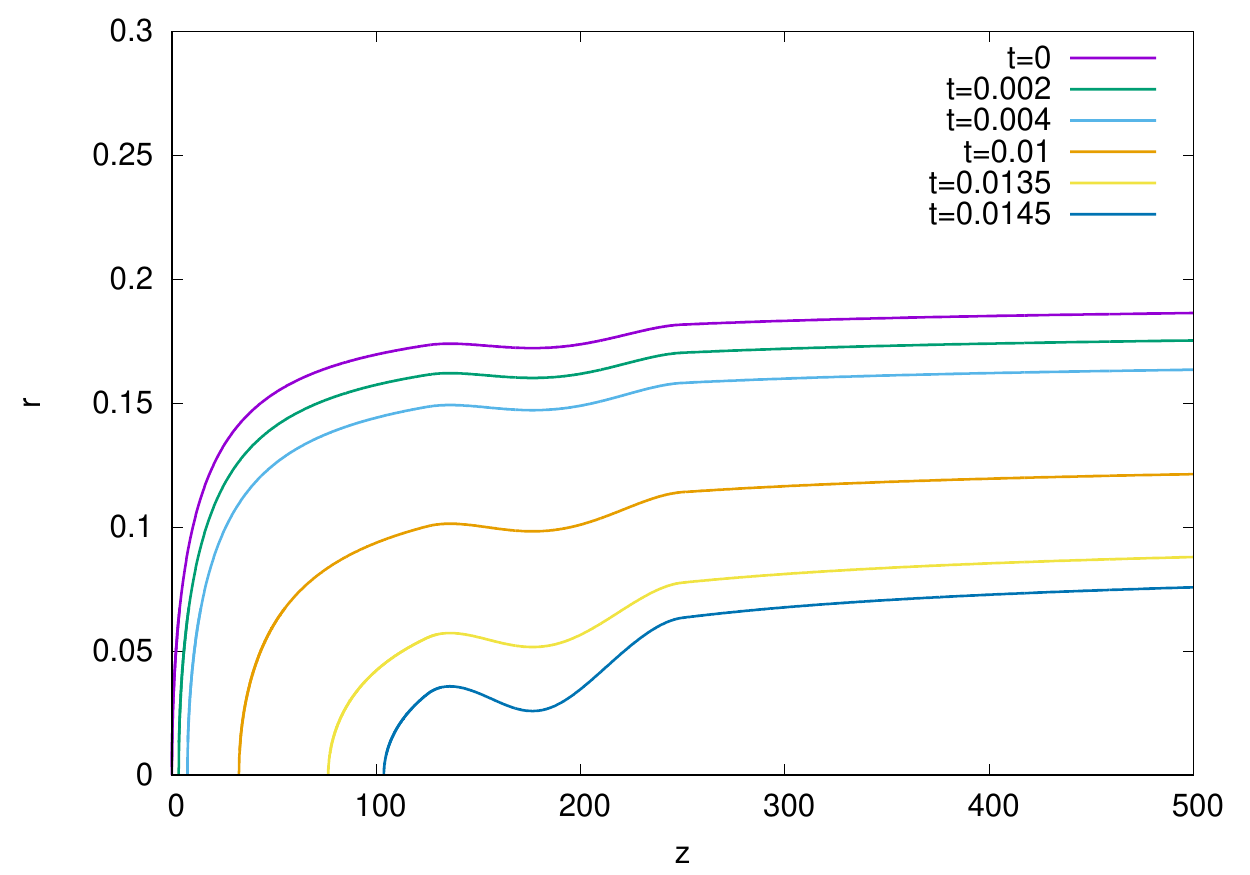}
	\caption{Numerical simulation of a MCF solution in the far class, for $n=3$  and $\gamma=3/4$.}
\end{figure}

\subsection{The critical class}

The numerical results for the near class solutions (cf. Section \ref{near_class}) and those for the far class solutions (cf. Section \ref{far_class}) display distinct singular behaviors. If we interpolate between the perturbations that lead to the near class and those that lead to the far class, then we expect to encounter critical behavior, corresponding to a closed interval of intermediate parameter values $[s_0, s_1]$ for $s_0\leq s_1$. The existence of such critical behavior is already observed and confirmed for the degenerate neckpinch in MCF of closed rotationally symmetric hypersurfaces \cite{AAG, AV97}. Indeed, if we cinch a round sphere at the equator, then we find that loose cinching leads to a global Type-I round singularity modelled by the sphere, whereas tight cinching leads to a local Type-I nondegenerate neckpinch modelled by the cylinder. Interpolating between these two classes of Type-I singular behaviors leads to a degenerate neckpinch exhibiting both Type-I and Type-II singular behaviors. However, to the best of the authors' knowledge, whether or not degenerate neckpinches in compact MCF occur for a closed interval of parameters or just a single threshold value is still unknown. We note that a similar phenomenon appears in the rotationally symmetric Ricci flows on spheres \cite{GI05, GI08, GZ08, AIK11, AIK15}.

\begin{figure}[H]
	\includegraphics[width=0.6\textwidth]{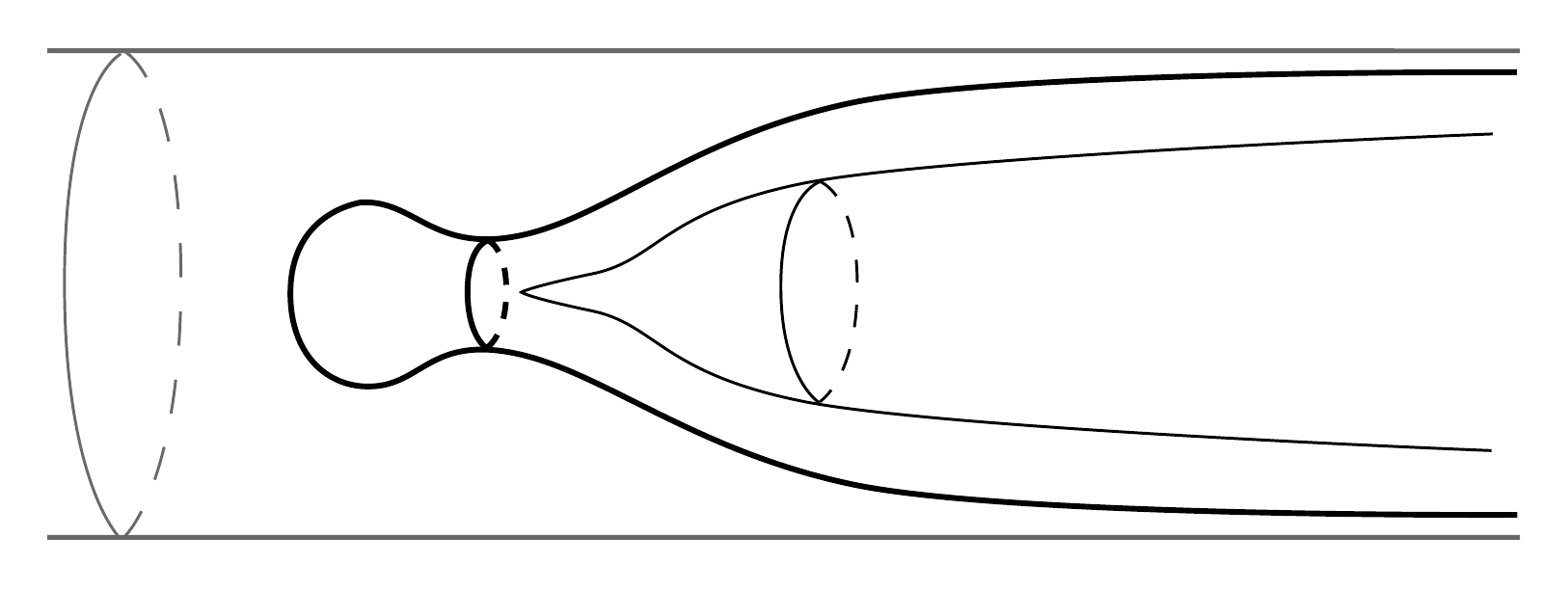}
	\caption{The conjectured behavior for a MCF solution in the critical class.}\label{critical_class}
\end{figure}

We conjecture the following behavior for a noncompact MCF solution in the critical class; see Figure \ref{critical_class} for a schematic drawing. If we perturb an initial hypersurface to have a mildly pinched neck at some critical distance from the tip, then we expect that the MCF of the region containing the tip and the neck should shrink to a cusp in the same fashion that MCF of an asymmetrical dumb-bell develops a degenerate neckpinch (cf. \cite{AV97}). The curvature near the neck is conjectured to blow up at the Type-I rate $(T-t)^{-1/2}$ and the singularity is expected to be modelled by a cylinder. We conjecture that the curvature at the tip blows up at a Type-II rate $(T-t)^{-(1-1/m)}$, where $m\geq 3$ is an integer. We note in particular that the Type-II blowup rates $(T-t)^{-(1-1/m)}$ interpolate between the Type-I rate $(T-t)^{-1/2}$ and the Type-II rates $(T-t)^{-(\gamma+1/2)}$, where $\gamma\geq 1/2$, for the solutions in the near class and those constructed in \cite{IW19, IWZ20}.

\section{Conclusions and discussions}\label{discussion}

Using a geometrically natural and analytically novel overlap method, we have carried out a numerical investigation of the local stability of rotationally symmetric, complete noncompact MCF solutions with Type-II curvature blowup as constructed in \cite{IW19, IWZ20}. Our numerical results indicate two distinct singular behaviors for two distinct families of local perturbations: local stability of the Type-II curvature blowup (at the tip) for solutions in the near class solutions, and local instability that leads to Type-I nondegenerate neckpinch for solutions in the far class.

Our numerical results strongly suggest the existence of a critical class as one interpolates between the near class and the far class. Detection of such a critical class by direct numerical testing can be very difficult, in part because it is more difficult to resolve a singularity developing in the overlapping region for our overlap method (cf. Section \ref{num_method}). Nevertheless, by drawing an analogy from the setting of MCF of topological spheres, in which a degenerate neckpinch occurs as a critical class in interpolating through Type-I singular solutions \cite{AAG, AV97}, we pose the following conjecture about the geometric and analytical nature of a class of critical solutions.

\begin{conj}\label{conj_critical}
	A solution in the critical class develops a local Type-II degenerate neckpinch whose precise asymptotics are obtained in \cite{AV97}.
\end{conj}

Assuming the above conjecture is true, it is then intriguing to ask whether this mean curvature flow might be able to continue through the degenerate neckpinch singularity. We note that similar questions for Ricci flow have been studied. Angenent, Caputo and Knopf \cite{ACK12} construct smooth forward Ricci flow solutions of singular initial metrics resulting from rotationally symmetric nondegenerate neckpinches on $S^{n+1}$ (occurring at the singular time $T$). In particular, the norm of the Riemann curvature tensor $|\Rm|$ of these solutions decreases at a rate (precisely,  $\frac{\log|t-T|}{t-T}$ for $t\geq T$), which is slightly faster than the Type-I rate before the singularity. In comparison, Carson has constructed smooth Ricci flows emerging from rotationally symmetric degenerate neckpinches on $S^{n+1}$, and the curvature of these solutions decreases at the same rate at which it blows up \cite{Car16} (see also \cite{Car18}).

So far we have only considered perturbations, and hence the perturbed MCF solutions, which are all rotationally symmetric. A natural next step for these numerical studies is to consider perturbations with a less restrictive assumption placed on the symmetry of these perturbations. In particular, we will be interested in non-rotationally symmetric perturbations near the tip where the curvature blow-up is Type-I. Then the numerical simulations will allow us to explore whether or not MCF of non-rotationally symmetric geometries evolve toward rotationally symmetric ones. We remark that for MCF of noncompact surfaces, the stability of a rotationally symmetric Type-I nondegenerate neckpinch under arbitrary $C^3$ perturbation is proved in a series of papers by one of the authors of this paper and his collaborators \cite{GK15, GKS18}. More generally, the stability of generalized cylinders in the class of self-shrinkers is proved by Colding and Minicozzi \cite{CM12}.

In this paper, we have investigated numerically the stability of rotationally symmetric, asymptotically cylindrical MCF solutions with Type-II curvature blowup in finite time. Generalizing to other asymptotic geometries, two of the authors of this paper and their collaborator have constructed MCF of entire graphs with super-linear growths at spatial infinity and a Type-II curvature blowup occurring in infinite time \cite{IWZ-2b}. A numerical stability analysis for these solutions can be pursued in the future.


\bibliography{mcf_type2_numerical}

\end{document}